\documentclass{siamltex}
\usepackage{graphicx}
\usepackage{amsmath,amssymb,leftidx}
\usepackage{bm}
\usepackage{bbm}
\usepackage[latin1]{inputenc}
\usepackage{todonotes}
\usepackage{algorithm,algorithmic}
\usepackage{color}
\usepackage{subfigure}
\usepackage{url}  
\usepackage{mathtools}

\DeclareMathOperator{\range}{Range}   %

\newcommand{\cA}{\mathcal{A}}
\newcommand{\cB}{\mathcal{B}}
\newcommand{\cC}{\mathcal{C}}

\newcommand{\cF}{\mathcal{F}}

\newcommand{\cH}{\mathcal{H}}

\newcommand{\cS}{\mathcal{S}}

\newcommand{\RR}{\mathbb{R}}

\newcommand{\Gr}{\mathrm{Gr}^3}

\renewcommand{\a}{\alpha}
\renewcommand{\b}{\beta}
\newcommand{\g}{\gamma}

\newcommand{\<}{\langle}
\renewcommand{\>}{\rangle}

\newcommand{\x}{\times}

\newcommand{\tp}{^{\sf T}}

\newcommand{\tml}[3][]{\bm{\left(} #2 \bm{\right)}_{ #1} \bm{\cdot} #3}
\newcommand{\tmr}[3][]{#2 \bm{\cdot} \bm{\left(} #3 \bm{\right)}_{ #1}}

\newcounter{excntr}[section]

\newenvironment{example}[1][\bigskip \textsc{\small EXAMPLE}]{\begin{trivlist} \refstepcounter{excntr}
\item[\hskip \labelsep  {\itshape #1} \arabic{section}.\arabic{excntr}]}{\endproof\end{trivlist}}

\title{Spectral Partitioning of Large and  Sparse Tensors using Low-Rank Tensor
  Approximation} 
\author{Lars Eld\'en\thanks{Department of Mathematics, Link\"{o}ping
    University, SE-581 83, Link\"{o}ping, Sweden
    ({lars.elden@liu.se})} \and Maryam Dehghan\thanks{Department of Mathematics, Persian Gulf University, 75169, Bushehr, Iran ({ma.dehghan@mehr.pgu.ac.ir, maryamdehghan880@yahoo.com})}\\
   }

\begin{document}

\maketitle

   \begin{center}
   \small  \today
   \end{center}

   \bigskip
   
\begin{abstract}
  The problem of partitioning  a large and sparse tensor is considered,
  where the tensor consists  of a sequence of adjacency
  matrices. Theory is developed that  is a generalization of 
  spectral graph partitioning. A best rank-$(2,2,\lambda)$
  approximation is computed for $\lambda=1,2,3$, and the partitioning
  is computed from the orthogonal matrices and the core tensor of the
  approximation. It is shown that if the tensor has a certain
  reducibility structure, then the solution of the best approximation
  problem  exhibits the reducibility structure of
  the tensor.  Further, if the tensor is close to being
  reducible, then still the solution of the exhibits the structure of
  the tensor. Numerical examples with synthetic data corroborate the
  theoretical results. Experiments with tensors from applications show
  that the method can be used to extract  relevant information from
  large, sparse, and noisy data.
\end{abstract}



\section{Introduction}
\label{sec:intro}

Spectral graph partitioning is a standard method in data science graph
applications, see, e.g., \cite{chung97,lux07,maca09,urzi14}. It is
usually based on the computation of the two smallest eigenvalues and
corresponding eigenvectors of  the graph Laplacian or its normalized
version. In view of the abundance in data science applications of
graph data organized as tensors, it is  natural to ask whether it is
possible to generalize spectral graph partitioning to sequences of
graphs. To our knowledge there is no natural generalization
of the graph Laplacian to tensors. However, the normalized Laplacian
approach is equivalent to a method based on the computation of the two
largest eigenvalues and corresponding eigenvectors of the normalized
adjacency matrix of the graph. If the graph is close to being
disconnected, then this property shows in the structure of the eigenvectors,
and can be  used to partition the graph. This method can be generalized to
tensors.   

Let $\cA$ be a tensor that consists of a sequence of normalized
adjacency matrices for undirected graphs; thus the matrices are
symmetric. Assume that almost all the graphs are close to being
disconnected in the same way (we will make these connectivity
properties   precise in Section \ref{sec:reducibility}). 
The proposed method is based on the computation of the best 
rank-$(2,2,\lambda)$ approximation
\begin{equation}
  \label{eq:best-appr-0}
  \min_{\mathrm{rank}(\cC)=(2,2,\lambda)} \| \cA - \cC \|,
\end{equation}
for $\lambda=1,2,3$, (for detailed definitions see Section
\ref{sec:concepts}). Thus we are 
considering low-rank Tucker 
approximations \cite{hit:27,tu:64} with best approximation
properties \cite{lmv:00b}. Let the  solution be  
\[
  \cC = \tml{U,U,W}{\cF},
\]
where $U\tp U = I_2$ and $W\tp W = I_\lambda$, and $\cF$ is a core
tensor of dimension $(2,2,\lambda)$. The connectivity properties of
$\cA$ can be deduced from the structure of $U$, $W$, and $\cF$, which
can be used for partitioning the tensor.

The best approximation problem \eqref{eq:best-appr-0} can be written,
equivalently, as the solution of
the maximization problem
\begin{equation}
  \label{eq:rayleigh}
   \max_{\Gr(2,2,\lambda)} \| \tmr{\cA}{X,X,Z} \|,
\end{equation}
where $\Gr(2,2,\lambda)$ is a product of certain Grassmann
manifolds. The maximization is a generalization of the  Rayleigh quotient
maximization of a symmetric matrix.

In Section \ref{sec:concepts} we  define notation and
preliminaries. To motivate the generalization to tensors, we give a
brief introduction to spectral graph partitioning using the normalized
adjacency matrix in Section \ref{sec:undir-graph}. A graph is
disconnected if the adjacency matrix is reducible. The corresponding
tensor concept is introduced in Section \ref{sec:spectral-tensor},
where we also recall the concept of stationary points for the
maximization problem \eqref{eq:rayleigh}. In Section
\ref{sec:best-reduc} we give  characterizations of the
solutions of the best approximation of tensors with certain
reducibility structures. These are used in a method for
computing the reducibility structure of tensors that are close to
reducible. We give examples of such synthetic tensors and compute their
structure in Section \ref{sec:examples}. There we also show how the
method can be used to analyse a real data tensor, which is far from
being reducible. 

The theory in this paper is based on the perturbation theory in
\cite{elsa11}, which shows  that if the tensor is perturbed by a small
amount then the solution $U$, $W$, and $\cF$ of
\eqref{eq:best-appr-0} are also perturbed by a small amount. Therefore
it is essential to use the best rank-$(2,2,\lambda)$ approximation.

The present paper is the second in a series of three on the analysis
of large and sparse tensor data. In the first \cite{eldehg20b} we
develop a Krylov-Schur like algorithm for efficiently computing best
rank-$(2,2,\lambda)$ approximations of large and sparse tensors. In
the third paper \cite{eldehg20c} we use the methods to analyze large
tensors from data science applications.

Other approaches to tensor clustering are described in 
\cite{jsb09,cwhyl13,wzzy19,cxz20,xzxgg20}, and in \cite{ldjd10}
clustering is done using 
HOSVD (Higher Order Singular Value Decomposition). 

\section{Tensor Concepts and Preliminaries}
\label{sec:concepts}
We will use the term tensor in a restricted sense, i.e. as a
generalization of the matrix concept:  
a tensor $\cA$ is 
3-dimensional array of real numbers, $\cA \in \RR^{l \times m \times
  n}$, where the vector space of tensors is equipped with  algebraic
structures to be defined.

Tensors will be denoted by calligraphic letters, e.g
$\mathcal{A},\mathcal{B}$, matrices by capital roman letters and
vectors by lower case roman letters.  In order not to burden the
presentation with too much detail, we sometimes will not explicitly
mention the dimensions of matrices and tensors, and assume that they
are such that the operations are well-defined. The whole presentation
will be in terms of tensors of order three, or equivalently 3-tensors.

 The different ``dimensions'' of the tensor
are referred to as \emph{modes}.  We will use both standard subscripts
and MATLAB-like notation: a particular tensor element will be
denoted in two equivalent ways, $\cA(i,j,k) = a_{ijk}$.

A subtensor obtained by fixing one of the indices is called a
\emph{slice}, e.g., $\cA(i,:,:)$.
In the case when the index is fixed in the first mode, we call the slice 
a $1-$slice, and correspondingly for the other modes.  
A  slice  can be  considered as an order-3 tensor
(3-tensor) with a singleton mode, and also as a matrix.  
A \emph{fiber} is a subtensor, where all indices but one are
fixed. For instance, $\cA(i,:,k)$ denotes a mode-2 fiber. 

For a given third order tensor, there are three associated subspaces,
one for each mode. These subspaces are given by  
\begin{align*}
& \range \{ \cA(:,j,k) \; | \; j = 1:m, \; k = 1:n  \},\\
& \range \{ \cA(i,:,k) \; | \; i = 1:l, \; k = 1:n  \},\\
& \range \{ \cA(i,j,:) \; | \; i = 1:l, \; j = 1:m  \}.
\end{align*}
The \emph{multilinear rank} \cite{hit:27,sili08} of the tensor is said
to be equal to $(p,q,r)$ 
if the dimension of these subspaces are $p,$ $q,$ and $r$, respectively. 

It is customary in numerical linear algebra to write out column
vectors with the elements arranged vertically, and row vectors with
the elements horizontally. This becomes inconvenient when we are
dealing with more than two modes. Therefore we will not make a
notational distinction between mode-1, mode-2, and mode-3 vectors, and
we will allow ourselves to write all vectors organized vertically. It
will be clear from the context  to which mode the vectors belong.
However, when dealing with matrices, we will often talk of 
them as consisting of column vectors.

Since we will be dealing with 3-tensors only, it is still possible to
display them graphically in figures and formulas. For instance, we
will consider tensors,
\[
  \RR^{2 \times 2 \times n} \ni
  \begin{pmatrix}
    p_1 & p_3\\
    p_3 & p_2
  \end{pmatrix}, \qquad p_i \in \RR^n, \quad i=1,2,3.
\]
It may help the reader's intuition, if the mode-3 fibers are thought
of as extending perpendicular to the plane of the paper/screen.


\subsection{Tensor-Matrix Multiplication}
\label{sec:ten-matmult}

We define \emph{multilinear multiplication of a tensor by a matrix} as
follows.  For concreteness we first present multiplication by one matrix
along the first mode and later for all three modes simultaneously. The
mode-$1$ product of 
a tensor $\cA \in \RR^{l \times m \times n}$ by a matrix $U \in \RR^{p
  \times l}$ is defined\footnote{The notation \eqref{eq:contra}-\eqref{eq:mat-tensor} was
suggested in  \cite{sili08}. An alternative notation was
  earlier given in \cite{lmv:00a}.  Our $\tml[d]{X}{\cA}$ is the same
  as $\cA \times_d X$ in that system.} 
\begin{equation}\label{eq:contra}
 \RR^{p \times m \times n} \ni \mathcal{B} =  \tml[1]{U}{ \cA} ,\qquad
b_{ijk} = \sum_{\a=1}^{l} u_{i \a} a_{ \a jk}  .
\end{equation}
This means that all mode-$1$ fibers in the $3$-tensor $\cA$ are
multiplied by the matrix $U$. Similarly, mode-$2$ multiplication by a
matrix $V \in \RR^{q \x m}$ means that all mode-$2$ fibers are
multiplied by the matrix $V$.  Mode-$3$ multiplication is analogous.
With a third matrix $W\in \RR^{r \x n}$, the tensor-matrix
multiplication in all modes is given by
\begin{equation} 
  \label{eq:mat-tensor}
\RR^{p \x q \x r}	\ni \cB = \tml{U,V,W}{\cA}, \qquad b_{ijk} = \sum_{\a,\b,\g=1}^{l,m,n} u_{i \a} v_{j \b} w_{k \g}a_{ \a \b \g},
\end{equation}
where the mode of each multiplication is understood from the order in
which the matrices are given.

It is convenient to introduce a separate notation for multiplication
by a transposed matrix $\bar{U} \in \RR^{l \times p}$:
\begin{equation}
  \label{eq:cov}
\RR^{p \times m \times n} \ni \mathcal{C}= \tml[1]{\bar{U}\tp }{\cA}
=\tmr[1]{\cA}{\bar{U}}, \qquad
c_{ijk} = \sum_{\a=1}^{l} a_{\a jk} \bar{u}_{\a i}.
\end{equation}
Let $u \in \RR^l $ be a  vector and $\cA \in \RR^{l \times m
  \times n}$ a tensor. Then
\begin{equation}
  \label{eq:ident-1}
 \RR^{1 \times m   \times n} \ni \cB := \tml[1]{u\tp}{\cA} =
\tmr[1]{\cA}{u} \equiv B \in \RR^{m   \times  n}.
\end{equation}
Thus we identify a tensor with a singleton dimension with  a
matrix. Similarly, with $u \in \RR^l$ and $w \in \RR^n$, we will
identify
\begin{equation}
  \label{eq:ident-13}
\RR^{1 \times m \times 1} \ni \cC :=   \tmr[1,3]{\cA}{u,w}  \equiv
c \in \RR^m,
\end{equation}
i.e., a tensor of order three with two singleton dimensions is
identified with a vector, here in the second mode.

\subsection{Inner Product, Norm, and Contractions}
\label{sec:norm}

Given two tensors $\cA$ and $\cB$ of the same dimensions, we define
the \emph{inner product},
\begin{equation}\label{eq:inner-prod}
\< \cA , \cB \> = \sum_{\a,\b,\g} a_{\a \b \g} b_{\a \b \g}.
\end{equation}
The  corresponding \emph{tensor norm} is
\begin{equation}
  \label{eq:norm}
  \| \cA \| = \< \cA, \cA \>^{1/2}.
\end{equation}
This \emph{Frobenius norm} will be used throughout the paper.
As in the matrix case, the norm is invariant under orthogonal
transformations, i.e.,
\[
  \| \cA \| = \left\| \tml{U,V,W}{\cA} \right\| = \| \tmr{\cA}{P,Q,S} \|,
\]
for orthogonal matrices $U$, $V$, $W$, $P$, $Q$, and $S$. This is obvious 
from the fact that multilinear multiplication by orthogonal matrices
does not change the Euclidean length of the corresponding fibers of the tensor.

For convenience we will denote the inner product of
  vectors $x$ and $y$ in any mode (but, of course, the same) by $x\tp
  y$.

The following well-known result will be needed.

\begin{lemma}\label{lem:LS}
  Let $\cA \in \RR^{l \times m \times n}$ be given along with three
  matrices with orthonormal columns, $U \in \RR^{l \times p}$, $V \in
  \RR^{m \times q}$, and $W \in \RR^{n \times r}$, where $p \leq l$,
  $q \leq m$, and $r \leq n$. Then the least squares problem
\[
\min_{\cS} \| \cA - \tml{U,V,W}{\cS} \|,
\]
with $\cS \in \RR^{p \times q \times r},$ has the unique solution
\[
\cS= \tml{U\tp ,V\tp ,W\tp }{\cA} = \tmr{\cA}{U,V,W}.
\]
The elements of the tensor $\cS$ are given by 
\begin{equation}
  \label{eq:elem-S}
  s_{\lambda\mu\nu}=\tmr{\cA}{u_\lambda,v_\mu,w_\nu}, \quad 1 \leq
  \lambda \leq p, \quad 1 \leq \mu \leq q, \quad 1 \leq \nu \leq r.
\end{equation}
\end{lemma}
  The proof is a straightforward generalization of the corresponding
  proof for matrices, see \cite{lmv:00b}.

The inner product \eqref{eq:inner-prod} can be considered as a special
case of the \emph{contracted product of two tensors},
which is a tensor (outer) product followed by a contraction along specified
modes. Thus, if $\cA$ and $\cB$ are $3$-tensors,
we define, using essentially the notation of \cite{bako:06a},

\addtocounter{equation}{1}
\begin{align}
  \cC &= \left\< \cA, \cB \right\>_{1}\,, &
   c_{jkj'k'} &= \sum_\a a_{\a j k} b_{\a j' k'}\,, & &
   \mbox{($4$-tensor)}\,, \label{eq:cntr1}\tag{\theequation.a}\\
  D &= \left\< \cA, \cB \right\>_{1,2}\,, &
   d_{k k'} & = \sum_{\a,\b} a_{\a \b k} b_{\a \b k'}\,, & &
   \mbox{($2$-tensor)}, \label{eq:cntr2}\tag{\theequation.b}\\
  e & = \< \cA, \cB \> = \left\< \cA, \cB \right\>_{1,2,3}\,, & e &
   = \sum_{\a,\b, \g} a_{\a \b \g} b_{\a \b \g}\,, & &
   \mbox{(scalar)}. \label{eq:cntr3}\tag{\theequation.c}
\end{align}
It is required that the contracted dimensions in the two
tensors are equal. We will refer to the first two as \emph{partial
  contractions}. It is also
convenient to introduce a notation when contraction is performed in all but
one mode. For example the product in \eqref{eq:cntr2} may also be written  
\begin{equation}
\left\< \cA, \cB \right\>_{1,2} \equiv \left\< \cA, \cB \right\>_{-3}\,.
\end{equation}
Two tensors $\cA$ and $\cB$ of equal dimensions are said to be
\emph{$(i,j)-$orthogonal}, if 
\[
\< \cA , \cB\>_{i,j} = 0.
\]
In the case of  3-tensors, if two tensors  $\cA$ and $\cB$ are
$(1,2)$-orthogonal, we also say that they are $-3$-orthogonal, i.e. 
$\< \cA, \cB \>_{-3} = 0.$

The definition of contracted products  is valid also  when the 
tensors are of different order. The only assumption is that the dimension
of the correspondingly contracted  modes are the same in the two
arguments. The dimensions of the resulting product are in the order
given by the non-contracted modes of the first argument followed by
the non-contracted modes of the second argument.

\subsection{Best Rank-$(r_1,r_2,r_3)$  Approximation}
\label{sec:best-appr}
The problem of approximating  tensor $\cA \in \RR^{l \x m \x n}$ by
another tensor $\cB$ of lower multi-linear rank, 
\begin{equation}
  \label{eq:best-appr-min}
\min_{\rank(\cB)=(r_1,r_2,r_3)} \| \cA - \cC \|^2,  
\end{equation}
is treated in \cite{lmv:00b,zhgo:00,elsa09,idav08}. It is shown in
\cite{lmv:00b} that \eqref{eq:best-appr-min} is equivalent to 
\begin{equation}
  \label{eq:best-appr-max}
  \max_{X,Y,Z} \| \tmr{\cA}{X,Y,Z} \|^2, \quad \mbox{subject to} \quad
  X\tp X = I_{r_1}, \quad
  Y\tp Y =I_{r_2}, \quad Z\tp Z = I_{r_3},
\end{equation}
where $X \in \RR^{l \x r_1}$, $Y \in \RR^{m \x r_2}$, and 
$Z \in \RR^{n \x r_3}$.  The problem \eqref{eq:best-appr-max} can be
considered as a \emph{Rayleigh quotient maximization problem}, in
analogy with the matrix case \cite{amsd:02}. We will refer to the two equivalent problems
as the \emph{the best rank-$(r_1,r_2,r_3)$ approximation problem}.
To
simplify the terminology somewhat we will refer to a solution
$(U,V,W)$ of the approximation problem as \emph{the best
  rank-$(r_1,r_2,r_3)$ approximation of the tensor}, and often
tacitly assume that the corresponding core tensor is $\cF =
\tmr{\cA}{U,V,W}$.  

The constrained maximization problem \eqref{eq:best-appr-max} can be
thought of as an unconstrained maximization problem on 
\emph{a product of Grassmann manifolds}, cf. \cite{eas98,ams:07}. 
 Thus a   solution $(U,V,W)$
  is a representative of  an equivalence class of 
  matrices $(U Q_1, V Q_2, W Q_3)$, where $Q_i \in \RR^{r_i \times r_i}$
  are orthogonal.  Note also that 
  \begin{equation}
    \label{eq:equiv}
\tml{U,V,W}{\cF} =
  \tml{U_\#,V_\#,W_\#}{\cF_\#},    
  \end{equation}
  where $(U_\#,V_\#,W_\#) =(U Q_1, V Q_2, W Q_3)$, and $\cF_\# =
  \tmr{\cF}{Q_1,Q_2,Q_3}$. For simplicity,  in the sequel when we talk
  about \emph{the solution} of the best approximation problem, it is 
   understood as \emph{a representative of the solution equivalence
    class}. Likewise, a  \emph{unique solution} means that the
  equivalence class is unique.

  It is also important to keep in mind that as the optimization
  problem is nonconvex, we cannot guarantee that a local optimum is a
  global one, see \cite{frge20}.

\subsection{(1,2)-Symmetric Tensors}
\label{sec:symmetric}

 A 3-tensor $\cA \in \RR^{m \x m \x n}$ is  called 
\emph{(1,2)-sym\-met\-ric} if all its 3-slices are symmetric, i.e.,
\[
\cA(i,j,k) = \cA(j,i,k), \qquad i,j=1,2,\ldots,m, \quad k=1,2, \ldots, n. 
\]
Henceforth in this paper we will assume that $\cA$ is {(1,2)-sym\-met\-ric}.  


We will partition  the (1,2)-symmetric  tensor with respect to the first two
modes, 
\[
\cA =
\begin{pmatrix}
  \cA_{11} & \cA_{12}\\
\cA^{'}_{12} & \cA_{22}
\end{pmatrix},
\]
where $\cA_{11} \in \RR^{m_1 \times m_1 \times n}$ for $m_1<m$, and the
  other ``block tensors'' have dimensions accordingly. Here we use the
  notation $\cA^{'}_{12}$ to indicate that each 3-slice is the
  transpose of the corresponding slice of $\cA_{12}$. 

  \medskip
  
\begin{proposition}\label{prop:symm-best}
  Assume that the tensor $\cA$ is (1,2)-symmetric, and that  the
  solution $(U,V,W)$ 
  of the  best approximation problem \eqref{eq:best-appr-min} with
  $r_1=r_2$ is unique. Then representatives of the solution can be
  chosen such that  $U=V$ and  the core tensor 
  $\cF=\tmr{\cA}{U,U,W}$  is   (1,2)-symmetric.  
\end{proposition}

\begin{proof}
    $\cA$ is (1,2)-symmetric, and, to be optimal, so is the
  approximating tensor $\tml{U,V,W}{\cF}$. Therefore $U$ and $V$  must
  span the same subspace, and we can choose $V=UQ$, for an orthogonal matrix
  $Q$. The rest follows from the Grassmann property of the problem,
  and from Lemma \ref{lem:LS}.  
\end{proof}
\medskip

In this paper we will treat the best approximation problem 
\begin{equation}
  \label{eq:best-appr-symm}
  \max_{X,Z} \| \tmr{\cA}{X,X,Z} \|, \qquad X\tp X = I_{r_1}, \quad
   Z\tp Z = I_{r_3},
\end{equation}
and denote its solution $(U,U,W)$. We  use the notation
\begin{equation}
  \label{eq:grassmann-sym}
  \Gr(r_1,r_1,r_3) = \{ (X,X,Z) \; | \; X\tp X = I_{r_1}, \; Z\tp Z =
  I_{r_3} \},
\end{equation}
in connection with (1,2)-symmetric minimization over a product of Grassmann
manifolds\footnote{The standard notation would be $\mathrm{Gr}(r_1,m) \times
  \mathrm{Gr}(r_1,m) \times \mathrm{Gr}(r_3,n)$. As the dimensions of
  the vector spaces ($m$ and $n$) 
  are implicit in the problem formulation, we use
  \eqref{eq:grassmann-sym},  for simplicity.}.
%

\subsection{A Krylov-Schur Like Method for Computing Best
  Rank-$(r_1,r_2,r_3)$    Approximations}
\label{sec:krylov-schur}
To be able to solve large problems with sparse tensors it is necessary
to use methods that do not apply transformations to the tensor itself
(except possibly reorderings), because transformations would
generally create fill-in. For problems with sparse matrices,
restarted Krylov methods are standard \cite{lesoya:98,stew:02}. In
\cite{sael11,eldehg20b} we developed a block-Krylov type method for
tensors, which accesses the tensor only in tensor-matrix
multiplications, where the matrix consists of a relatively small
number of columns. Thus the method is memory-efficient.  In
\cite{eldehg20b} it was combined with the Krylov-Schur approach for
the computation of the best rank-$(r_1,r_2,r_3)$ approximation of a large
and sparse tensors.
For the problems solved in
Section \ref{sec:examples} convergence is very fast. It is shown in
\cite{eldehg20b} that the block Krylov-Schur method is in general
considerably faster than the 
Higher Order Orthogonal Iteration \cite{lmv:00b}.

\section{Background: Spectral Graph Partitioning Using the  Adjacency   Matrix}  
\label{sec:undir-graph}

In this section we give a brief overview of the ideas behind spectral
graph partitioning, emphasizing the aspects that we will generalize to
tensors. A survey of methods is given in \cite{lux07}, and
the mathematical theory is presented in \cite{chung97}. 
 The  presentation here is based on the  fact that
graph partitioning using the normalized graph Laplacian is equivalent
to partitioning using the normalized adjacency matrix.


Let $A_0$ be the \emph{(unnormalized) adjacency matrix} of a connected
undirected graph $G$.   The \emph{Laplacian} of the graph is
\[
L_0 = D - A_0,
\]
where $D$ is a diagonal matrix, 
\[
D = \diag(d), \quad d = A_0 e,
\] 
and $e=(1 \, 1 \, \cdots 1)\tp$.   
The \emph{normalized Laplacian} is defined
\[
{L} = I - {A} : = I - D^{-1/2} A_0 D^{-1/2}.
\]
The matrix $A$ can be called the \emph{normalized adjacency matrix},
and it is symmetric. 
Both Laplacian matrices are positive semidefinite. 
Usually spectral partitioning is developed in terms of the Laplacian
or the normalized Laplacian \cite{fied73,chung97}.  The smallest
eigenvalue of either Laplacian is equal to zero, and the value of
the second smallest determines how close the graph is to being
disconnected. Obviously, the normalized Laplacian has the same
eigenvectors as the adjacency matrix $A$ and the eigenvalues are
related by 
\[
\lambda(L) = 1 - \lambda(A).
\]
So instead of basing spectral partitioning on the smallest eigenvalues
of $L$ (or ${L}_0$) we can consider the largest eigenvalues of ${A}$
and the corresponding eigenvectors. We will assume the ordering of the
eigenvalues of the normalized adjacency matrix $A$, 
\[
1=\lambda_1 \geq \lambda_2 \geq \cdots ,
\]
and denote the corresponding  eigenvectors $v_1, v_2, \ldots$.
 
\begin{figure}[htbp!]    
\centering
\includegraphics[width=1.0\textwidth]{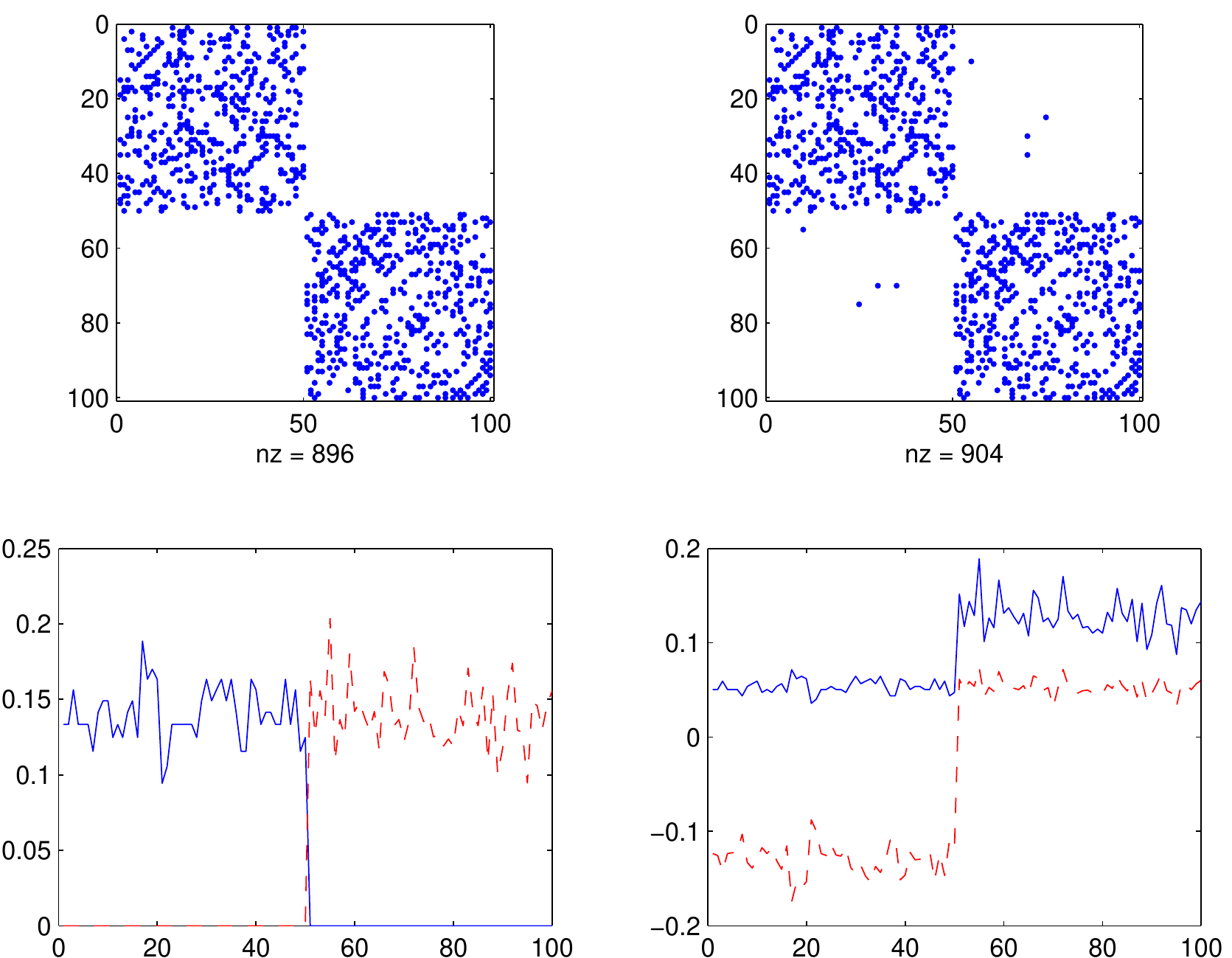}
\caption{Spy plots of the normalized adjacency matrices $A$ of a graph with two
  components (upper left), and $\widehat{A}$ of the same graph where
  four edges have
  been added between the two components (upper right). The first two
  eigenvectors of the adjacency matrices are given in the lower plots.}
\label{fig:AA0}
\end{figure}

To illustrate the heuristics of spectral graph partitioning, we
constructed a disconnected graph with 100 nodes, and a nearby graph,
where four edges were added that made the graph connected.  The
respective normalized adjacency matrices $A$ and $\widehat{A}$ are
illustrated in Figure \ref{fig:AA0}. Note that the matrix $A$ is
\emph{reducible} (actually in reduced form), whereas $\widehat{A}$ is
\emph{irreducible} (see, e.g., \cite[Section 6.2.21]{hojo:85}).

The following fundamental properties of the eigenvalue problem for
symmetric matrices are the key to spectral partitioning. We give the
results in somewhat loose form,  for precise statements, see, e.g.,
\cite[Theorems 8.1.4, 8.1.10]{govl13}, \cite[Theorems 17.2,
17.3]{eldenmatrix19}.

\begin{proposition}\label{prop:perturbation}
  If a symmetric matrix $A$ is perturbed by the amount $\epsilon$,
  then the eigenvalues are  perturbed by approximately the same
  amount  $\epsilon$. The eigenvectors are perturbed essentially by $\epsilon$
  divided by a quantity that depends on the distance to nearby
  eigenvalues. 
\end{proposition}

The adjacency matrices are normalized and therefore 
 the two-subgraph matrix $A$  has a double eigenvalue equal to
1, which are the largest eigenvalues for the adjacency matrices of each
sub-graph. The two largest eigenvalues of $\widehat{A}$ corresponding
to   the connected graph are 1
and $0.9913$.

The first  two eigenvectors of the two adjacency matrices are given in
Figure \ref{fig:AA0}. Since $A$ has a double eigenvalue, the
eigenvectors are not unique, only the eigenspace.  One can
choose the basis vectors for the eigenspace as the eigenvectors for
each of the two submatrices, which are nonnegative. 
For $\widehat{A}$   the first eigenvector is positive, due to the fact
that the matrix is irreducible and nonnegative (see Perron-Frobenius
theory, e.g. in \cite[Chapter8]{mey:00}). The
second eigenvector, which is orthogonal to the first, must have a sign
change, and that sign change occurs at the interface between the
diagonal blocks. Since $\widehat{A}$ is a small perturbation of $A$, the
subspaces spanned by the first two eigenvectors are close; here the
cosines of the angles between the subspaces are $1.0000$ and $0.9997$.

In real applications, given a graph that is close to being
disconnected, the nodes are not ordered so that the structure is
apparent as in our example. In fact, the purpose of spectral
partitioning is to find the node ordering that  exhibits such
structure. A simplified formulation of  one step of spectral
partitioning of a connected  graph is  given in Algorithm
\ref{alg:graph-part}. 
\begin{algorithm}
\caption{One step of spectral graph partitioning}
\label{alg:graph-part}
\begin{algorithmic}
\STATE \textbf{1.}  Compute the first two
  eigenvectors $v_1$  and $v_2$ of the normalized  adjacency matrix. 

\STATE \textbf{2.} Order the elements of $v_2$ in decreasing order, and
reorder the   nodes correspondingly.  

\STATE \textbf{3.} Given that ordering, compute the cost  for cuts in the
  vicinity of the sign change of the reordered eigenvector, and
  partition the graph in two where the  cost is the smallest. 
\end{algorithmic}
\end{algorithm}

``Cost'' in the algorithm can be given a precise meaning, e.g., in
terms of the concept \emph{conductance}, see
\cite{chung97}. Essentially it should be a measure of the number of
edges that are broken at a cut.

Finally in this section we summarize a few facts about spectral graph
partitioning that will be the basis for generalization to the tensor
case. 

\begin{description}

\item{\textbf{A1.}} Using the normalized adjacency matrix is
  equivalent to using the normalized Laplacian. 

\item{\textbf{A2.}} If the graph is connected, then the first
  eigenvector $v_1$ is positive, and the second has 
  positive and negative components (due to orthogonality). 

\item{\textbf{A3.}} The matrix $B = \lambda_1 v_1 v_1\tp + \lambda_2 v_2 v_2\tp$  is
  the best rank-2 approximation of the normalized adjacency matrix. 

\item{\textbf{A4.}} In the disconnected case the normalization makes
  the two blocks each have the largest eigenvalue equal to 1, with
  nonnegative eigenvectors. 
\end{description}

\medskip

\section{Spectral Partitioning of Tensors}
\label{sec:spectral-tensor}

 We first  present two motivating examples, which    generalize
 the graph concept. The objects of study are sparse (1,2)-symmetric
 tensors.
 
\medskip

\begin{example}\label{ex:web-link}
  \textbf{Web link analysis.} In \cite{koba:06} a graph representing
  links between web pages is studied. Thus the first two modes
  represent the web pages.  Instead of just registering the existence
  of   links, a third mode is used to denote the key word (term)  via which
  the two web pages are linked, in the style of Wikipedia. 
Thus,  a 3-tensor $\cA$ is constructed as
\[
a_{ijk}=\begin{cases} 1 & \mbox{if page $i$  points to page $j$  using
    term $k$,} \\
       0 & \mbox{otherwise.}
        \end{cases}
\]
Assuming that the links are bidirectional,  each 3-slice of $\cA$
represent the undirected graph of pages that are 
  linked using one particular term. 

  {Assuming that there are  two  different topics covered on
  two (almost) disjoint sets of web pages, can we find those topics
  and the corresponding groups of web pages? Can we also find the key
  words that are common to all web pages?}  
\end{example}

\medskip
 
\begin{example}\label{ex:network-time}
  \textbf{ Network  evolving over  time.}  In  \cite{bhk07}
  a subset of the Enron e-mail communication is analyzed over a period
  of 44 months. The   first two tensor modes represent 184 employees 
  communicating, and the third mode is time. Thus each 3-slice is
  equivalent to  a   graph of the e-mail users who communicate during one
  time unit. 

  Assume that there are two groups of users, and that both groups
  communicate frequently internally during one period in time.  Then
  at a later stage in time, all users communicate with all the
  others\footnote{In \cite{dka11} such temporal dynamics is
    discovered using the time factors of a CP decomposition of the
    tensor. }.

  On the other hand, assume that there are two groups of users that
  both communicate only within their own group during the the  whole
  period studied.

  As a third example, assume that during the first time period
   one group communicates with all others, while the other group does
   not communicate within its own group. Then during the
   second period, all users are communicating with all.
   

  {Can we detect such communication patterns by spectral
    partitioning?} 
\end{example}

\medskip

In both examples, since we assume that the graphs are undirected, the
tensor is $(1,2)$-symmetric. Also, in both examples the nodes have no
loops, i.e., for all $i$ and $k$, $a_{iik}=0.$

\subsection{Reducibility}
\label{sec:reducibility}

Reducibility and Perron-Frobenius theory for tensors are  treated
in \cite{cpz08,fgh13,gahe16}. From 
Example \ref{ex:network-time}  one can see that we need a more specific 
reducibility concept than that in \cite{cpz08}. Thus we make the
following definition for  the case of (1,2)-symmetric 3-tensors. 

\begin{definition}\label{def:reducibility}
  A (1,2)-symmetric tensor $\cA \in \RR^{m \times
    m \times n}$   is called
  (1,2)-reducible if there   exists a nonempty proper subset $I \subset
  \{1,2,\ldots,m\}$  and a non\-empty subset $K \subseteq
  \{1,2,\ldots,n\}$ such that  
\[
a_{i j k} = a_{jik} = 0, \quad \forall i \in I,  \quad \forall j
\notin I, \quad \forall k \in K.
\]
Similarly, it is called 3-reducible if 
\[
a_{ijk} = 0, \quad \forall i,j \in I, \quad \forall k \in K. 
\]
If  $\cA$  is either (1,2)-reducible or 3-reducible we call it reducible,
and if it is not reducible, we  call it {irreducible}. 
\end{definition}

If the tensor is reducible, then we can permute its indices in a
(1,2)-symmetric way\footnote{I.e., so that (1,2)-symmetry is
  preserved.} so that $I$ becomes $\{1,2,\ldots,n_1\}$, for some 
$n_1 < n$, and also permute the indices in the third mode in a usually
different way so that $K$ becomes $\{1,2,\ldots,m_1\}$ for some $m_1
\leq m$. A (1,2)-reducible tensor will after such a permutation have
the structure
\begin{equation}
  \label{eq:reducible-tensor12}
  \cA(:,:,1:m_1) = 
  \begin{pmatrix}
    \cA_1 & 0 \\
      0   & \cA_2
  \end{pmatrix},
\quad 
  \cA(:,:,m_1:m) = 
  \begin{pmatrix}
    \cA_3 & \cA_4 \\
      \cA'_4   & \cA_5
  \end{pmatrix}.
\end{equation}
Note that we allow $m_1=m$, in which case the second subtensor in
\eqref{eq:reducible-tensor12} is vacuous.  A 3-reducible tensor has
the structure 
\begin{equation}
  \label{eq:reducible-tensor3}
  \cA(:,:,1:m_1) = 
  \begin{pmatrix}
    \cA_1 & \cA_2 \\
      \cA'_2& 0 
  \end{pmatrix},
\quad 
  \cA(:,:,m_1:m) = 
  \begin{pmatrix}
    \cA_3 & \cA_4 \\
      \cA'_4   & \cA_5
  \end{pmatrix},
\end{equation}
for some $m_1 < m$. (We here use the prime notation to indicate that
$\cA_4'(:,:,k)\tp =\cA_4(:,:,k)$ for $k=m_1+1,\ldots,m$). Note that
\eqref{eq:reducible-tensor12} corresponds to the first case in Example
\ref{ex:network-time} and \eqref{eq:reducible-tensor3} corresponds to
the third case.


\subsection{Normalization of 3-Slices}
\label{sec:normalization}
In Section \ref{sec:undir-graph} we described the normalization of an
adjacency matrix for a connected, undirected graph. There the normalization
has three main effects: (i) the largest eigenvalue of the adjacency
matrix is equal to 1, (ii) if the graph is close to being disconnected,
then the second largest eigenvalue is close to 1, (iii) the eigenvectors
corresponding to the two largest eigenvalues give information about
the structure of the graph, and can be used to find a good
partitioning.

We are not aware of any analogous normalization of a  (1,2)-symmetric
3-tensor with analogous  properties. However, if we assume that   the
3-slices are adjacency matrices of undirected graphs, it is reasonable
to normalize them separately. Thus we assume that the tensor has been
normalized as in Algorithm \ref{alg:normalize}. The vector $e=(1 \, 1
\, \cdots 1)\tp \in \RR^m$. 
\begin{algorithm}
\caption{Normalization of the 3-slices of (1,2)-symmetric $\cA$ }
\label{alg:normalize}
\begin{algorithmic}
\STATE \textbf{for} $i=1:n$

\STATE $\qquad B:= \cA(:,:,i)$

\STATE $\qquad d:=Be$ 

\STATE $\qquad D:=\diag(d)$

\STATE $\qquad \cA(:,:,i):=D^{-1/2}BD^{-1/2}$

\STATE \textbf{end} 
\end{algorithmic}
\end{algorithm}

\medskip

The 3-slices of a normalized (1,2)-symmetric tensor all have the
maximum eigenvalue equal to 1.



\subsection{Stationary Points}
\label{sec:stationary}

 For later reference we here restate the first order conditions for a
 stationary point for the maximization problem
 \eqref{eq:best-appr-symm} \cite[Proposition 4.2]{elsa11}.

Consider the the solution $(U,U,W)$ of the (1,2)-symmetric best
rank-$(r_1,r_1,r_3)$  approximation problem (i.e., $U \in \RR^{m \x r_1}$ and $W
\in \RR^{n \x r_3}$), and let $U_\perp$ and $W_\perp$ be matrices such
that $(U \; U_\perp)$ and $(W \; W_\perp)$ are both square and
orthogonal. 

\begin{proposition}\label{prop:first-order}
Let $(U,U,W)$ be a 
stationary point   for the (1,2)-symmetric best approximation
problem, and put $\cF = \tmr{\cA}{U,U,W}$. Then the following
equations are satisfied
 \begin{align}
  \< \cF_\perp^1, \cF \>_{-1} &= 0, \qquad \cF^1_\perp =
\tmr{\cA}{U_\perp,U,W},   \label{eq:F-orth-1}\\
\< \cF_\perp^2, \cF \>_{-2} & = 0, \qquad \cF^2_\perp =
\tmr{\cA}{U,U_\perp,W},   \label{eq:F-orth-2}\\
 \< \cF_\perp^3, \cF \>_{-3} &= 0, \qquad \cF^3_\perp =
\tmr{\cA}{U,U,W_\perp}.  \label{eq:F-orth-3}
 \end{align}
%
\end{proposition}

In our case, when $\cA$ is (1,2)-symmetric, the two relations
\eqref{eq:F-orth-1} and \eqref{eq:F-orth-2} are equivalent. 

If we apply Proposition \ref{prop:first-order} to the best
rank-$(1,1,1)$  approximation problem, we
have the following result 
\cite{zhgo:00}.

\begin{lemma}\label{lem:111F}
A stationary point $(u,v,w)$    for the best rank$-(1,1,1)$
approximation of a (non-symmetric) tensor $\cA$ satisfies the equations
\begin{equation}
  \label{eq:111nonsym}
  \tmr[1,2]{\cA}{u,v}= \tau w, \qquad \tmr[1,3]{\cA}{u,w}=\tau v,
  \qquad \tmr[2,3]{\cA}{v,w} = \tau u,
\end{equation}
where $\tau = \tmr{\cA}{u,v,w}$.

In the case when $\cA$ is (1,2)-symmetric the solution $(u,u,w)$
satisfies
\begin{equation}
  \label{eq:111sym}
  \tmr[1,2]{\cA}{u,u} = \tau w, \qquad \tmr[1,3]{\cA}{u,w}=
  \tmr[2,3]{\cA}{u,w} = \tau u. 
\end{equation}
\end{lemma}

\medskip

  The lemma follows from the orthogonality relations in
  Proposition \ref{prop:first-order}, and  the fact that $\cF$ is 
   scalar: then  $\RR^{(m-1) \times 1 \times 1} \ni \cF_\perp^1 =0$, etc.

In the rest of the paper we will use the fact  that if the tensor $\cA$ is
nonnegative, then the solution vectors in the best $\rank-(1,1,1)$
approximation problem are nonnegative. This  seems to be
well-known; for completeness we give a simple proof in  Appendix
\ref{app:111}.


\section{Best Approximation Problems for Reducible Tensors}
\label{sec:best-reduc}

We will now present the solution of best low-rank approximation
problems for a few specific (1,2)-symmetric reducible tensors, and we
will show how the reducibility structure of the tensors corresponds to
a certain structure of the solution.  We  remind the reader that
$(U,U,W)$ should be understood as an equivalence class represented by
$(U,U,W)$, cf. 
the arguments around \eqref{eq:equiv}.  In Section
\ref{sec:find-struct} we will consider tensors that are close to being
reducible, and show how the theory of this section can be used to
compute partitionings.

Assume that   the (1,2)-symmetric tensor $\cA$ is
(1,2)-reducible with 
\begin{equation}
    \label{eq:A1A2}
\cA =
\begin{pmatrix}
  \cA_{1} & 0 \\
     0     &   \cA_{2}
\end{pmatrix}, \qquad \cA_{\nu} \in \RR^{m_\nu \x m_\nu \x n}, \quad
\nu=1,2, \quad m_1 + m_2 = m.     
  \end{equation}
Thus we  consider the special case when $n_1=n$  in
\eqref{eq:reducible-tensor12}.

\begin{lemma}\label{lem:U-structure}
    Assume that the (1,2)-symmetric tensor $\cA$ has the structure
    \eqref{eq:A1A2}, where the subtensors $\cA_1$ and $\cA_2$  are
    irreducible, and  assume further  that $(U,U,W)$ is the unique solution 
    of the best rank-(2,2,2) approximation problem. Then $U$ has one
    of the structures
\begin{equation}\label{eq:U-struct1}
\begin{pmatrix}
  u_{11} & u_{12} \\
  0 &  0
\end{pmatrix}, \qquad 
\begin{pmatrix}
    0 & 0 \\
  u_{21} & u_{22}
\end{pmatrix}, 
\end{equation}
or
\begin{equation}
  \label{eq:U-structure}
\begin{pmatrix}
  u_{11} & 0 \\
  0 & u_{22}
\end{pmatrix},
\end{equation}
where  $u_{1i} \in \RR^{m_1}$, $u_{2i} \in \RR^{m_2}$, for $i=1,2$.

The same result holds for the best rank-(2,2,1) approximation of $\cA$. 
  \end{lemma}
  \begin{proof} 
    Let $(X,X,Z) \in \Gr(2,2,2)$,
    and let
    $\cH \in \RR^{2 \times 2 \times 2}$ be arbitrary
    $(1,2)$-sym\-met\-ric. Put $\cC = \tml{X,X,Z}{\cH}$, and partition
      \[
    \cC =
    \begin{pmatrix}
      \cC_1 & \cC_{12}\\
      \cC_{12}^{'} & \cC_2      
    \end{pmatrix}
  \]
  Clearly $\cC$ is $(1,2)-$symmetric.  
  Then
  %
  \begin{align}
   \min_{\Gr(2,2,2)} \| \cA - \cC \|^2 &= \min_{\Gr(2,2,2)} \left(\|
                                         \cA_1 - \cC_1 \|^2 + 
    2 \|  \cC_{12} \|^2 + \| \cA_2 - \cC_2 \|^2\right) \nonumber \\
    & \geq 
    \min_{\Gr(2,2,2)} \left(\| \cA_1 - \cC_1 \|^2 +     \| \cA_2 - \cC_2
      \|^2\right)       \label{eq:AB-resid2} \\
    &= \min_{\Gr(2,2,2)} \left\|
      \begin{pmatrix}
        \cA_1 & 0 \\
        0 & \cA_2
      \end{pmatrix} - 
    \begin{pmatrix}
      \cC_1 & 0 \\
      0 & \cC_2
    \end{pmatrix} \right\|^2,
  \end{align}
  with strict inequality if $\cC_{12} \neq 0$.  
   We can write the rank-(2,2,2) approximation
  \begin{equation*}
    \cC = \tml[3]{z_1}{H_1} +\tml[3]{z_2}{H_2},
  \end{equation*}
  where $H_\nu = \tml[1,2]{X,X}{\cH(:,:,\nu)}, \, \nu=1,2,$ are
  symmetric rank-2
  matrices.  Due to the orthogonality of $z_1$ and $z_2$, the
  (1,2)-block of $\cC$ is equal to zero if and only if the corresponding
  (1,2)-blocks of $H_1$ and $H_2$ are both zero.    From Lemma
  \ref{lem:rank-2} in Appendix \ref{app:C12} we then see that the
  solution $U$ must
  have the structure \eqref{eq:U-struct1} or \eqref{eq:U-structure}.

  The proof for the best rank-(2,2,1) approximation is almost identical. 
\end{proof}

\medskip

If, in the above proof, we let $\cC_1$ have rank $(2,2,2)$, then, due
to the structure of $\cC$, we must have $\cC_2=0$, and vice versa.  So
we can let $\cC_1$ have rank $(2,2,2)$, put $\cC_2=0$, and minimize
the first term in \eqref{eq:AB-resid2}. This corresponds to the first
structure in \eqref{eq:U-struct1}. Similarly, minimizing only the second
term in \eqref{eq:AB-resid2} corresponds to the second structure in
\eqref{eq:U-struct1}. Thirdly, we can let $\cC_1$ and $\cC_2$ each
have rank (1,1,2), and minimize both terms in \eqref{eq:AB-resid2}
simultaneously. This corresponds to the structure in
\eqref{eq:U-structure}.

\medskip

\begin{proposition}\label{prop:A1A2}
Let the  (1,2)$-$symmetric, nonnegative  tensor $\cA \in \RR^{m \times
  m \times n}$ have the structure  
\begin{equation}
    \label{eq:A1A2-1}
\cA =
\begin{pmatrix}
  \cA_{1} & 0 \\
     0     &   \cA_{2}
\end{pmatrix}, \qquad \cA_{i} \in \RR^{m_i \x m_i \x n}, \quad
i=1,2, \quad m_1 + m_2 = n,
  \end{equation}
  where $\cA_1$ and $\cA_2$ are irreducible.  Assume that
  \begin{equation}
    \label{eq:Condition}
   \min_{\Gr(2,2,2)} \| \cA - \tml{X,X,Z}{\cF} \|^2 <
    \min_{\Gr(2,2,2)} \| \cA_1 - \tml{X,X,W}{\cF} \|^2 + \| \cA_2 \|^2,
  \end{equation}
  and that the corresponding inequality holds with $\cA_1$ and $\cA_2$
  interchanged.  Further assume that 
  the solution  $(U,U,W)$ 
   with core tensor $\cF$, of the best rank-(2,2,2) app\-roximation problem for
  $\cA$,   is unique.  Then 
  \begin{equation}
    \label{eq:U-structure-prop}
    U=
    \begin{pmatrix}
      u_1 & 0 \\
      0   & u_2 
    \end{pmatrix}.
  \end{equation}
   The 3-slices of $\cF$ are
  diagonal.

  The corresponding result hold for the best rank-(2,2,1)
  approximation. 
\end{proposition}

\begin{proof}
  The condition \eqref{eq:Condition} and Lemma
  \ref{lem:U-structure}  rule out the
  possibility that the solution has the the structure
  \eqref{eq:U-struct1}. So it must have the structure
  \eqref{eq:U-structure-prop}. We then have 
  \begin{align*}
    \cF &= \tmr[3]{\left(\tmr[1,2]{\cA}{U,U}\right)}{W} =
    \tmr[3]{\begin{pmatrix}
      \tmr[1,2]{\cA_1}{u_1,u_1} & 0\\
      0 & \tmr[1,2]{\cA_2}{u_2,u_2}
    \end{pmatrix}}{W}\\
    &= \begin{pmatrix}
      \tmr{\cA_1}{u_1,u_1,W} & 0\\
      0 & \tmr{\cA_2}{u_2,u_2,W}
    \end{pmatrix} \in \RR^{2 \times 2 \times 2};
  \end{align*}
  thus the 3-slices of $\cF$ are diagonal.
  
  The proof for the rank-(2,2,1) case is almost identical.
\end{proof}

\medskip


The following corollary can be proved by a modification of the proof
of Proposition \ref{lem:nonneg} in Appendix \ref{app:111}.

\medskip

\begin{corollary}
  The solution vectors $u_1$ and  $u_2$   in Proposition
  \ref{prop:A1A2} are nonnegative. 
\end{corollary}

\medskip
We will refer to nonnegative  vectors with the structure
\eqref{eq:U-structure-prop} 
as vectors in \emph{indicator form}, as they exhibit  the reducibility
structure of the tensor.

If the 3-slices of the tensor $\cA$ are  normalized as  described in
Section \ref{sec:normalization}, then the two blocks $\cA_1$ and $\cA_2$
in Proposition \ref{prop:A1A2} have  approximately equal magnitude, in
the sense that all 3-slices of the blocks have largest eigenvalue
equal to 1. Therefore, if the blocks $\cA_1$ and $\cA_2$ are far from being
reducible, it is likely that the assumptions \eqref{eq:Condition} are
satisfied.

A  special case of that in Proposition \ref{prop:A1A2} has an explicit
solution. 

\medskip

\begin{proposition}\label{prop:A000000A}
  Let $\cA \in \RR^{m \times m \times n}$ be
  (1,2)-symmetric, nonnegative, and (1,2)-reducible, with the
  structure
  \begin{align*}
  &\cA(:,:,1:n_1) =
\begin{pmatrix}
  \cA_{1} & 0 \\
     0     &   0
   \end{pmatrix},
   \qquad
     \cA(:,:,n_1+1:n) =
\begin{pmatrix}
   0 & 0 \\
     0     &  \cA_2
\end{pmatrix},\\
 &\cA_{i} \in \RR^{m_i \x m_i \x n_i}, \quad
i=1,2, \qquad m_1 + m_2 = m, \quad n_1+n_2=n,
  \end{align*}
where $\cA_1$ and $\cA_2$ are irreducible.
Let $(u_i,u_i,w_i)$  be the solutions of the best
rank-(1,1,1) approximation problem for $\cA_i$, $\, i=1,2$,
respectively, and let $\tmr[1,2]{\cA_i}{u_i,u_i}=\tau_i w_i, \;
i=1,2$. Assume that
\begin{equation}
  \label{eq:cond-A0000}
  \max_{\Gr(2,2,2)} \| \tmr{\cA_i}{X,X,Z} \|^2 < \tau_1^2 + \tau_2^2,
  \quad i=1,2.
\end{equation}
Then
the solution of the best rank-(2,2,2) problem for $\cA$ is
\begin{equation*}
  U =
  \begin{pmatrix}
     u_1 & 0 \\
    0 & u_2 
  \end{pmatrix},
\end{equation*}
and 
\begin{equation}\label{eq:W-ind}
  W =
  \begin{pmatrix}
    w_1 & 0\\
    0      & w_2
  \end{pmatrix}.
\end{equation}
The core tensor is
\[
  \cF(:,:,1)=
  \begin{pmatrix}
    \tau_1 & 0\\
    0 & 0
  \end{pmatrix}, \quad
  \cF(:,:,2)=
  \begin{pmatrix}
    0 & 0\\
    0 & \tau_2
  \end{pmatrix}.
\]
\end{proposition}

The proof is easy and we omit it here. Note that $U$ and $W$ are
nonnegative. 
 
Next assume that the 3-slices of the tensor $\cA$ are adjacency matrices of
bipartite graphs. This a special case of 3-reducibility, and  the
tensor has the structure as in the following proposition. 

\begin{proposition}\label{prop:0BB0}
  Let  the  (1,2)$-$symmetric, nonnegative tensor $\cA \in \RR^{m
    \times m \times   n}$ have the structure
  \[
    \cA =
    \begin{pmatrix}
      0 & \cB \\
      \cB^{'} & 0
    \end{pmatrix}, \qquad \cB \in \RR^{m_1 \times m_2 \times n}, \quad
    m_1 + m_2 =m.
  \]
  Assume that the best rank-$(1,1,1)$ approximation of
  $\cB$ is unique, given by nonnegative
  $(u,v,w)$, and let $\tmr[1,2]{\cB}{u,v} = \tau
  w$. Define the matrix $U$,
  \begin{equation}
    \label{eq:U-struct}
    U =
    \begin{pmatrix}
      0 & u\\
       v & 0
    \end{pmatrix}.
  \end{equation}
  Then the best rank-$(2,2,1)$ approximation of $\cA$ is given by
  $(U,U,w)$, and the core tensor is 
  \[
    \RR^{2 \times 2 \times 1} \ni \cF =
    \begin{pmatrix}
      0 & \tau \\
      \tau & 0
    \end{pmatrix}.
  \]
  
\end{proposition}

\begin{proof}
  We will show that $(U,U,w)$ satisfies
  \eqref{eq:F-orth-1}-\eqref{eq:F-orth-3}, and thus is a stationary
  point.  
  Choose $\bar U_\perp$, $\bar V_\perp$, and $\bar W_\perp$ such that
  $(u \, \bar U_\perp)$, $(v \, \bar V_\perp)$, and $(w \, \bar W_\perp)$
  are orthogonal matrices. From \eqref{eq:111nonsym} we see that
  \[
    \cF^1_{1\perp} := \tmr{\cB}{\bar U_\perp,v,w} = 0, \quad
    \cF^2_{2\perp} := \tmr{\cB}{u,\bar V_\perp,w} = 0, \quad
    \cF^3_{3\perp} := \tmr{\cB}{u,v,\bar W_\perp} = 0.
  \]
  Clearly
  \[
    \RR^{(m-2) \times m} \ni U_\perp =
    \begin{pmatrix}
      0 & \bar U\perp\\
      \bar V_\perp& 0
    \end{pmatrix}
 \]
 satisfies $U\tp U_\perp=0$. 
  Therefore
  \begin{align*}
    \cF_\perp^1 & = \tmr{\cA}{U_\perp,U,w} =
                  \begin{pmatrix}
                    0 & \tmr{\cB}{u,\bar V_\perp,w}\\
                    \tmr{\cB}{\bar U_\perp,v,w} & 0
                  \end{pmatrix}
     =
      \begin{pmatrix}
        0 & \cF_{2\perp}^2\\
        \cF_{1\perp}^1 & 0
      \end{pmatrix}
      = 0,
  \end{align*}
  where we have used the identity
  $\tmr{\cB'}{\bar V_\perp,u,w}=\tmr{\cB}{u,\bar V_\perp,w}$.
  Similarly one can prove $\cF_\perp^3=0$.  It is easy to prove
  that the point is optimal: Let $\cC = \tml{X,X,z}\cH$, and consider
  \[
    \left\|
      \begin{pmatrix}
        0 & \cB\\
        \cB' & 0 
      \end{pmatrix}
      -
      \begin{pmatrix}
        \cC_{11} & \cC_{12}\\
        \cC_{12}' & \cC_{22}
      \end{pmatrix}
      \right\|^2 \geq 
  \left\|
      \begin{pmatrix}
        0 & \cB\\
        \cB' & 0 
      \end{pmatrix}
      -
      \begin{pmatrix}
        0 & \cC_{12}\\
        \cC_{12}' & 0
      \end{pmatrix}
    \right\|^2
    = 2 \| \cB - \cC_{12} \|^2.
  \]
  Here  $\rank(\cC_{12})$ must be equal to  $(1,1,1)$, because otherwise
  \[
      \begin{pmatrix}
        0 & \cC_{12}\\
        \cC_{12}' & 0
      \end{pmatrix}
    \]
    would have rank higher than (2,2,1). Thus the optimum is achieved
    when $\cC_{12}$ is the best rank-(1,1,1) approximation of
    $\cB$. The expression for 
  $\cF = \tmr{\cA}{U,U,w}$ follows directly from the identity
  $\tmr[1,2]{\cB}{u,v} = \tau w$.
\end{proof}

\medskip

We will  now give an heuristic argument why, under certain
assumptions,  the best 
rank-(2,2,2) approximation problem for a tensor that is
(1,2)-reducible in the general sense of Definition
\ref{def:reducibility}, has a solution, where $U$ and $W$ are close to
being in indicator form.    Consider  $\cA \in \RR^{m \times m \times n}$,
which is    (1,2)-symmetric, nonnegative, (1,2)-reducible, with the
structure 
  \begin{equation}
    \label{eq:A-1,2-red-1}
    \cA =
    \begin{pmatrix}
      \cA_1 & \cB \\
      \cB' & \cA_2
    \end{pmatrix}, \qquad \cA_1 \in \RR^{m_1 \times m_1 \times n},
    \quad \cA_2 \in \RR^{m_2 \times m_2 \times n},  
      \end{equation}
where $m_1 + m_2 =m$,    and 
\begin{equation}    \label{eq:A-1,2-red-2}
  \cA(:,:,1:n_1) 
=  \begin{pmatrix}
    \cA_1^{(1)} & 0 \\
      0   & \cA_2^{(1)}
  \end{pmatrix}, \qquad 
  \cA(:,:,n_1+1:n)) =   \begin{pmatrix}
    \cA_1^{(2)} & \cB^{(2)} \\
      {\cB^{(2)}}'   & \cA_2^{(2)}
  \end{pmatrix},
\end{equation}
for some $n_1$ satisfying  $1< n_1 < n$.
 Assume further that $\cA_1$, $\cA_2$,
 and $\cB^{(2)}$  are irreducible. If
 \[
   \| \cA_i \| >> \| \cB \|, \quad \mathrm{or} \quad \| \cB \| >> \|
   \cA_i \|, \quad i=1,2,
 \]
 then, from Propositions \ref{prop:A1A2} and \ref{prop:0BB0}, and based
 on the perturbation 
 arguments in Section \ref{sec:find-struct}, the best
 rank-(2,2,2) and rank-(2,2,1) approximations will have $U$ in
 app\-roximate  indicator form. So we will assume that $\| \cA_i \|
 \approx \| \cB  \|$ for $i=1,2$.
 
  Then, with obvious notation for the partitioning of the
  approximating rank-(2,2,2) tensor $\cC$,
  \begin{align}\label{eq:C-1}
    \| \cA - \cC \|^2 
    &= \left\|
      \begin{pmatrix}
        \cA_1^{(1)} - \cC_1^{(1)}  & - \cC_{12}^{(1)}\\
        (- \cC_{12}^{(1)})' & \cA_2^{(1)} - \cC_2^{(1)}
      \end{pmatrix}
                              \right\|^2
       + \left\|
      \begin{pmatrix}
        \cA_1^{(2)} - \cC_1^{(2)}  & \cB^{(2)} - \cC_{12}^{(2)}\\
        (\cB^{(2)}- \cC_{12}^{(2)})' & \cA_2^{(2)} - \cC_2^{(2)}
      \end{pmatrix}
                                 \right\|^2\\
       &= \| \cA_1 - \cC_1 \|^2 + 2 \| \cC_{12}^{(1)} \|^2
         + \| \cA_2 - \cC_2 \|^2 + 2 \| \cB^{(2)} - \cC_{12}^{(2)}
         \|^2
         \label{eq:C-2}\\
       &\geq \| \cA_1 - \cC_1 \|^2 
          + \| \cA_2 - \cC_2 \|^2 + 2 \| \cB^{(2)} - \cC_{12}^{(2)}
        \|^2.\label{eq:C-3}         
  \end{align}
  %
%
%
  We can  write   the rank-(2,2,2) approximation
  \[
    \cC = \sum_{\nu=1}^2
    \tml[3]{w_\nu}{\left(\tml[1,2]{U,U}{\cF(:,:,\nu)}\right)},
    \]
    where $\cF(:,:,\nu)$ can be considered as 2-by-2
    matrices. Clearly, the optimal solution should have $\|
    \cC_{12}^{(1)}\|$ small or equal to zero.   From
    Lemma \ref{lem:rank-2} we see that $\cC^{(1)}_{12}=0$ if
\begin{equation}\label{eq:U-struct-11}
  U  =
  \begin{pmatrix}
     u_1 & 0 \\
     0 & u_2 
   \end{pmatrix}, \qquad
   W=
   \begin{pmatrix}
     w_1 & 0 \\
     0   & w_2
   \end{pmatrix},
\end{equation}
and $\cF(:,:,1)$ is diagonal. With $U$ and $W$ of this structure  we
can write the objective function
    \begin{align*}
      \| \cA - \cC \|^2 &= \| \cA_1^{(1)} -
       \tml{u_1,u_1,w_1}{\cF(1,1,1)}\|^2 
       + \| \cA_2^{(1)} -
             \tml{u_2,u_2,w_1}{\cF(2,2,1)}\|^2\\ 
                 &+ \| \cA_1^{(2)} - 
                   \tml{u_1,u_1,w_2}{\cF(1,1,2)}\|^2  
               + \| \cA_2^{(2)} 
            -\tml{u_2,u_2,w_2}{\cF(2,2,2)}\|^2 \\ 
      &+ 2 \| \cB^{(2)} -
        \tml{u_1,u_2,w_2}{\cF(1,2,2)}\|^2.
    \end{align*}
    Thus each tensor block is approximated by a rank-(1,1,1)
    tensor, which can be expected to give a rather small minimum. On
    the other hand, it is possible that the minimum in the 
    coupled minimization of \eqref{eq:C-1} will
    be smaller if the value of $\| \cC_{12}^{(1)} \|$  is a little larger
    than zero. This happened in  numerical examples in Section
    \ref{sec:examples}. However, our numerical experiments indicate
    that in many cases the minimization forces  $\| \cC_{12}^{(1)} \|$
    to be so small that the solution will not be in exact indicator
    form \eqref{eq:U-struct-11}, but in \emph{approximate indicator
      form}, so that the reducibility structure of the tensor is
    clearly exhibited.

\section{Computing the Approximate Structure of    Tensors that are Close to Reducible}
\label{sec:find-struct}

It is safe to assume that in real  applications exactly reducible data
tensors occur very seldom. In this section we will discuss
how the results on reducible tensors can be used to ascertain the
structure of tensors that are close to being reducible. This is based
on the theory of Section \ref{sec:best-reduc} and a perturbation
result for the best rank-$(r_1,r_2,r_3)$ approximation for tensors
that is an analogue to Proposition \ref{prop:perturbation} for
matrices. To avoid burdening the presentation with technical details,
we give a rather loose formulation of the results.

\begin{proposition}\cite[Section 5]{elsa11}
  \label{prop:perturbation-tensor}
  If a tensor  $A$ is perturbed by $\epsilon$,
  then the core tensor $\cF$ of the best
rank-$(r_1,r_2,r_3)$ approximation is  perturbed by approximately the same
  amount  $\epsilon$. The solution matrices $(U,V,W)$  are
  perturbed essentially by $\epsilon$ 
  divided by a quantity that is similar to  the distance to nearby
  eigenvalues in the eigenvalue case (the ``gap''). 
\end{proposition}

It is shown in \cite{elsa11} that the perturbation theory for the
singular value decomposition and the symmetric eigenvalue problem are
special cases of  Proposition \ref{prop:perturbation-tensor}. 

Proposition \eqref{prop:perturbation-tensor} shows that, if a
(1,2)-symmetric tensor is close to reducible, then the qualitative
properties of the solution $(U,U,W)$ and $\cF$ of the best
rank-$(r_1,r_1,r_3)$ approximation are close to those of the exactly
reducible nearby tensor. Thus, $(U,U,W)$ and $\cF$ can be used to
investigate the structure of the tensor, and partition it accordingly.

In the discussion on structure determination we first consider tensors
that are exactly reducible. Clearly, we must handle the Grassmann
non-uniqueness of the solution of the best approximation problem.  Let
the best rank-$(r_1,r_1,r_3)$ approximation be
$\cC = \tml{U,U,W}{\cF}$. Then, for any orthogonal
$Q_1 \in \RR^{r_1 \times r_1}$ and $Q_3 \in \RR^{r_3 \times r_3}$,
\[
  \cC = \tml{UQ_1,UQ_1,WQ_3}{(\tmr{\cF}{Q_1,Q_1,Q_3})}.
\]
For the reducible problems we have considered  in Propositions
\ref{prop:A1A2}, \ref{prop:A000000A}, and \ref{prop:0BB0},
we have written the solution  in indicator form, 
\eqref{eq:U-structure} and \eqref{eq:W-ind}.   If we multiply a pair
of vectors in indicator form by a rotation matrix, we get
\begin{equation}
  \label{eq:U-rot}
  U Q_1 =
  \begin{pmatrix}
    u_1 & 0\\
    0 & u_2
  \end{pmatrix}
  \begin{pmatrix}
    c & s \\
    -s & c
  \end{pmatrix}
  =
  \begin{pmatrix}
    c u_1 & s u_1\\
    -s u_2 & c u_2
  \end{pmatrix}. 
\end{equation}
Thus, due to the Grassmann indeterminacy, we can expect the solution
to be in the form of the right hand side of \eqref{eq:U-rot}. 
Since $u_1 \geq 0$ and $u_2 \geq 0$, one of the columns remains
nonnegative (or nonpositive) after multiplication, while the other has
positive and 
negative elements (the same happens if we multiply by a rotation
matrix with determinant $-1$). This implies that we can read off the
structure of the tensor by looking at the positive and negative
entries of one of the vectors.  For the tensor in non-reduced form,
say $\tml{P_1,P_1,P_3}{\cA}$, for permutation matrices $P_1$ and
$P_3$, the rows of \eqref{eq:U-rot} are permuted by $P_1$; thus, the
structure of the tensor can easily be found by reordering in
increasing order the column of $UQ_1$  that has both negative and
positive entries, and reordering the tensor accordingly.

But the positive-negative pattern of one of the vectors of $U$ (and
also $W$), is not enough.  To determine if the tensor is reducible, we
need to transform the computed $U$ and $W$ to indicator form, and
transform the core $\cF$ tensor accordingly; if the core has one of
the structures in the propositions in Section \ref{sec:best-reduc},
then the tensor is reducible. Thus we want to determine the rotation
that restores the indicator form.

Assume that the rows of $U$ have been reordered so that the negative
elements of are  at the bottom of the second column and that the
first  column is nonnegative. We want to find the rotation that puts
the vectors in indicator form: 
\[
  \begin{pmatrix}
    u_{11} & u_{12}\\
    u_{21} & u_{22}
  \end{pmatrix}
  \begin{pmatrix}
    c & -s\\
    s & c
  \end{pmatrix}
   =
   \begin{pmatrix}
     * & 0 \\
     0 & *
   \end{pmatrix}.
 \]
 A simple computation gives
 \begin{equation}
   \label{eq:Uhat}
   \widehat{U} q
   :=
   \begin{pmatrix}
     u_{12} & - u_{11}\\
     u_{21} & u_{22}
   \end{pmatrix}
   \begin{pmatrix}
     c\\
     s
   \end{pmatrix}
   =
   \begin{pmatrix}
     0\\
     0
   \end{pmatrix}.
 \end{equation}
 The vector $q$ is the right singular vector corresponding to the
 smallest singular value  of $\widehat{U}$. Thus, we can transform $U$
 to indicator form if the smallest singular value of $\widehat U$ is
 equal to zero, and the transformation is determined from the
 corresponding right singular vector. 

In principle,  the same procedure can be used if the tensor is close
to reducible. Assume that we have determined orthogonal matrices $Q_1$ and $Q_3$
that transform the computed $U$  and $W$ to approximate indicator form
(in the rank-(2,2,1) case $Q_3=1$ and $W$ is not transformed). The 
corresponding tensor is $\tmr{\cF}{Q_1,Q_1,Q_3}$, and we will refer to
it as the \emph{structure tensor}.

An algorithm for computing the structure of  a (1,2)-symmetric
tensor can be based on the successive computation of best
rank-(2,2,1), rank-(2,2,2), and rank-(2,2,3) approximations, along
with indicator vectors and structure tensor, and using
the characterisations of the solutions for different structures in
Section \ref{sec:best-reduc}. This is illustrated in Section
\ref{sec:examples}.

\section{Examples}
\label{sec:examples}

We  first give  four  examples with synthetic data that illustrate the
reducibility structures in the Propositions in Section
\ref{sec:best-reduc}, and describe how our  algorithm 
finds the structure. Then we apply the method to a tensor with real
data and show how topics in news text  can be extracted using the partitioning.
Further examples with real data  are given in \cite{eldehg20c}. 

\subsection{Synthetic Data}
\label{sec:synth}

For the examples we created random, sparse,
non\-nega\-tive and 
(1,2)-symmetric tensors in $\RR^{200 \times 200 \times 200}$ with
certain reducibility structures, and added a (1,2)-symmetric
perturbation, sparser than the original tensor. The perturbation made
the tensor irreducible, but close to reducible. All tensor elements
were put equal to 1, before we normalized the 3-slices as described in
Section \ref{sec:normalization}. We then performed a (1,2)-symmetric
permutation, and, in some cases, another permutation in the third
mode. We refer to this as the permuted tensor. Then we restored the
reducibility pattern via a best low-rank approximation. The
approximation was computed by the Krylov-Schur-like BKS
algorithm \cite{eldehg20b} with
relative termination criterion $10^{-12}$. After having been
initialized with one HOOI iteration for rank-$(1,1,1)$  
\cite{lmv:00b},  the BKS  method converged in
3-7 iterations.

The experiments were performed using Matlab on a desktop computer with
8 Gbyte primary memory. To give a rough idea about the amount of work, we
mention that the execution time for Examples \ref{ex:Ex1}-\ref{ex:Ex4}
was between 1.5  and 3.5 seconds.

In order to investigate how robust the reconstruction of reducility
structure is, we made the perturbed tensors relatively far from
reducible. However, this makes it slightly more difficult to draw
conclusions from the structure tensor. On the other hand, the
indicator vectors always gave the correct structure. In all three
examples, when we made the perturbation smaller, the structure tensors
became closer to those in Section \ref{sec:best-reduc}.  With no
perturbation in Examples \ref{ex:Ex1} and \ref{ex:Ex2} the vectors
were in indicator form to working precision.

\medskip

\begin{example}\label{ex:Ex1}
  We constructed a tensor with the structure close to that  in  Proposition
  \ref{prop:A1A2}.  The original  tensor is illustrated in Figure
\ref{fig:OrigEx1}, where we have used a 3D equivalent of the  Matlab
\texttt{spy} plot. To exhibit the sparsity, we also give a spy plot of
a 3-slice. 
%
%
\begin{figure}[htbp!]    
  \centering
  \includegraphics[width=.5\textwidth]{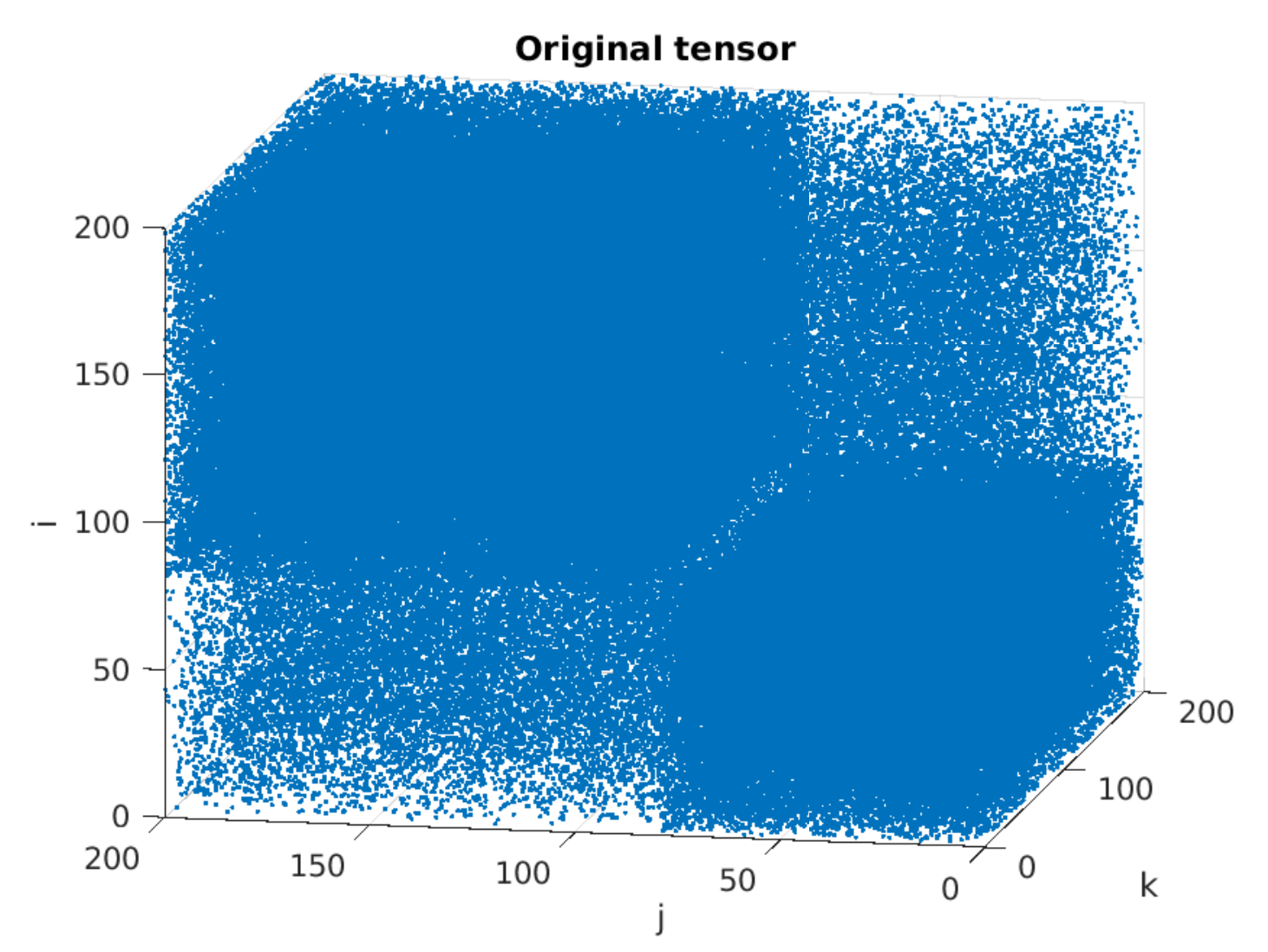}
\includegraphics[width=.35\textwidth]{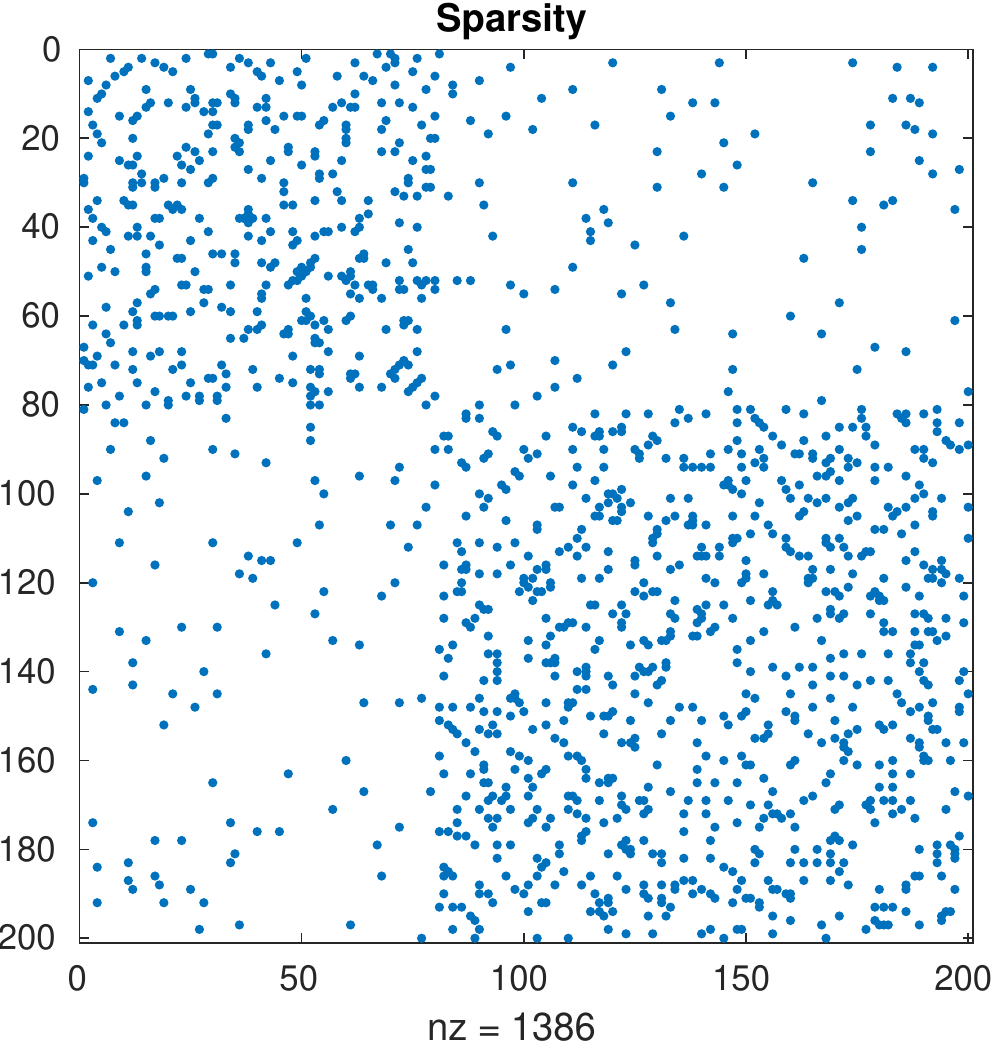}
\caption{Example \ref{ex:Ex1}. Original tensor (left),  a 3-slice 
  (right).}
\label{fig:OrigEx1}
\end{figure}
We  computed the best rank-(2,2,1) approximation, reordered $U(:,2)$
to make it decreasing,
computed indicator vectors and structure tensor. The vectors are
illustrated in Figure \ref{fig:VectorsEx1}. 
\begin{figure}[htbp!]    
\centering
\includegraphics[width=.65\textwidth]{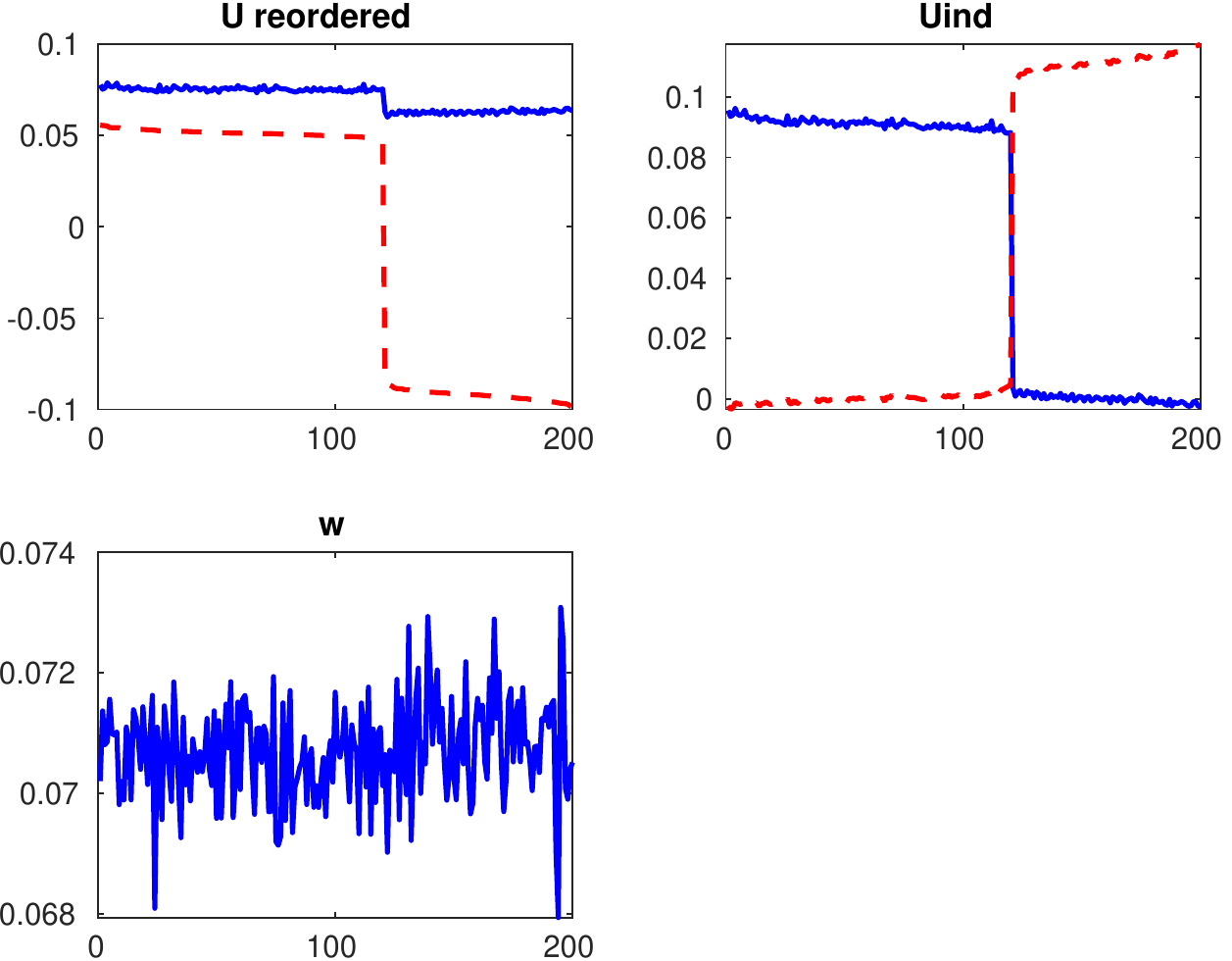}
\caption{Example \ref{ex:Ex1}, rank-(2,2,1) approximation. The column
  vectors of the solution 
   $U$ after reordering that makes $U(:,2)$ decreasing (upper left), 
   indicator 
  vectors (upper right), and the vector $w$  (lower left). $U(:,1)$ is
solid and blue, $U(:,2)$ is dashed and red. } 
\label{fig:VectorsEx1}
\end{figure}

Applying the reordering of the rows of $U$ to the permuted tensor, we
got  the one 
illustrated in Figure \ref{fig:ReorderedEx1}. 
\begin{figure}[htbp!]      
\centering
\includegraphics[width=.6\textwidth]{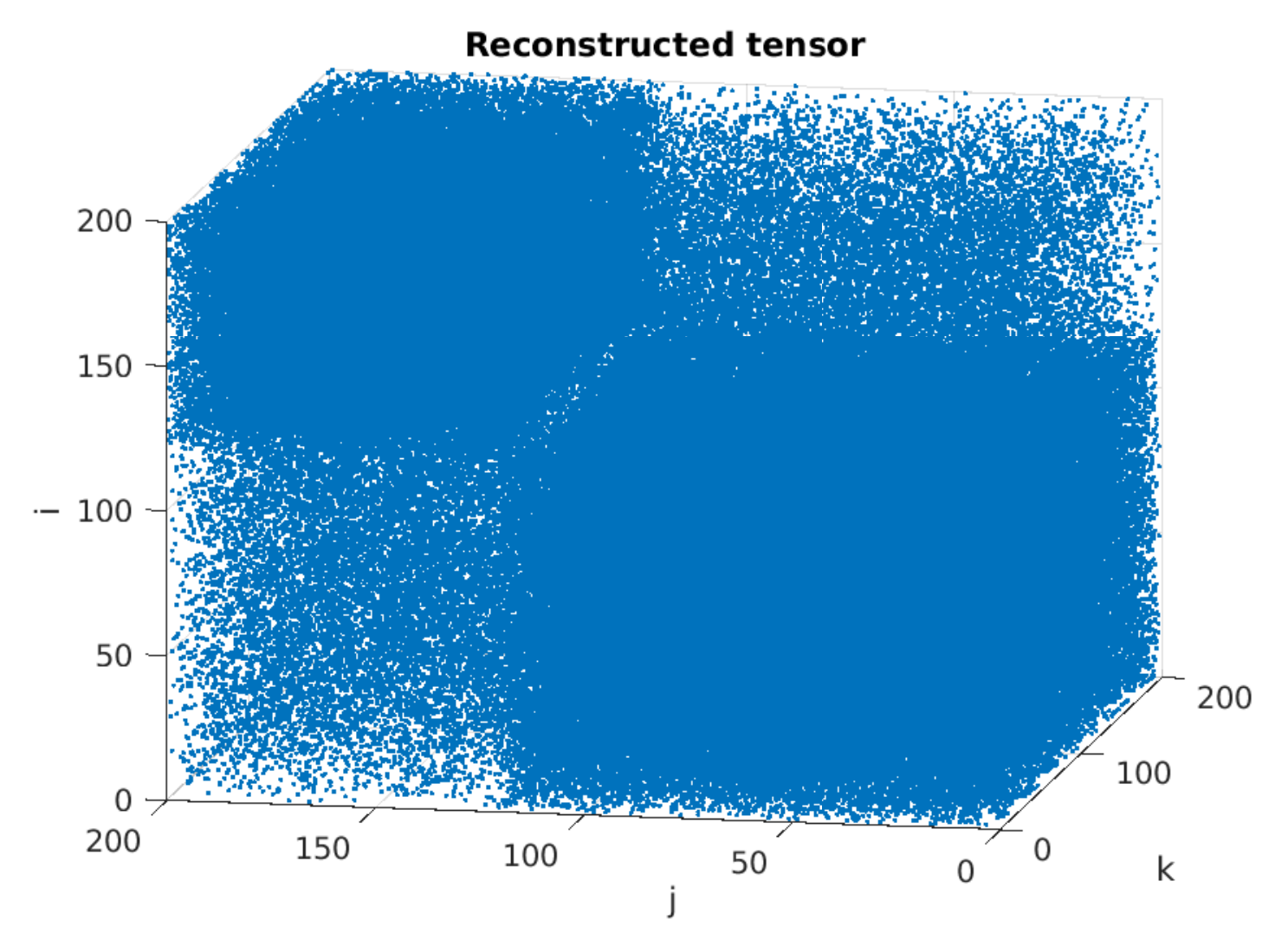}
\caption{Example \ref{ex:Ex1}. The reconstructed tensor.}
\label{fig:ReorderedEx1}
\end{figure}
The structure is seen already in the reordered vectors, but before we can draw
any safe conclusion must look at  the
structure tensor,
\[
    \cF_{\mathrm{struct}} =
    \begin{pmatrix}
          12.29 &  -1.95\\ 
	   -1.95 &  10.75
    \end{pmatrix}.
  \]
  This indicates that we have identified the correct structure.  The
  norm of  the core tensor and the original tensor were $\| \cF \| =
  16.56$ and $\| \cA \| = 76.21$.  To further check
  if there was any undetected significant structure in the third mode we also
  computed the best rank-(2,2,2) approximation. The norm of the core
  tensor now was 16.57, which shows that practically no structure
  was present that required a rank-(2,2,2) approximation. Moreover,
  while the $U$ vectors remained similar to 
  those in Figure \ref{fig:VectorsEx1}, the computed indicator vectors
  in Figure \ref{fig:VectorsEx1-222}
  for $W$ did not reveal any third mode structure.
\begin{figure}[htbp!]    
\centering
\includegraphics[width=.65\textwidth]{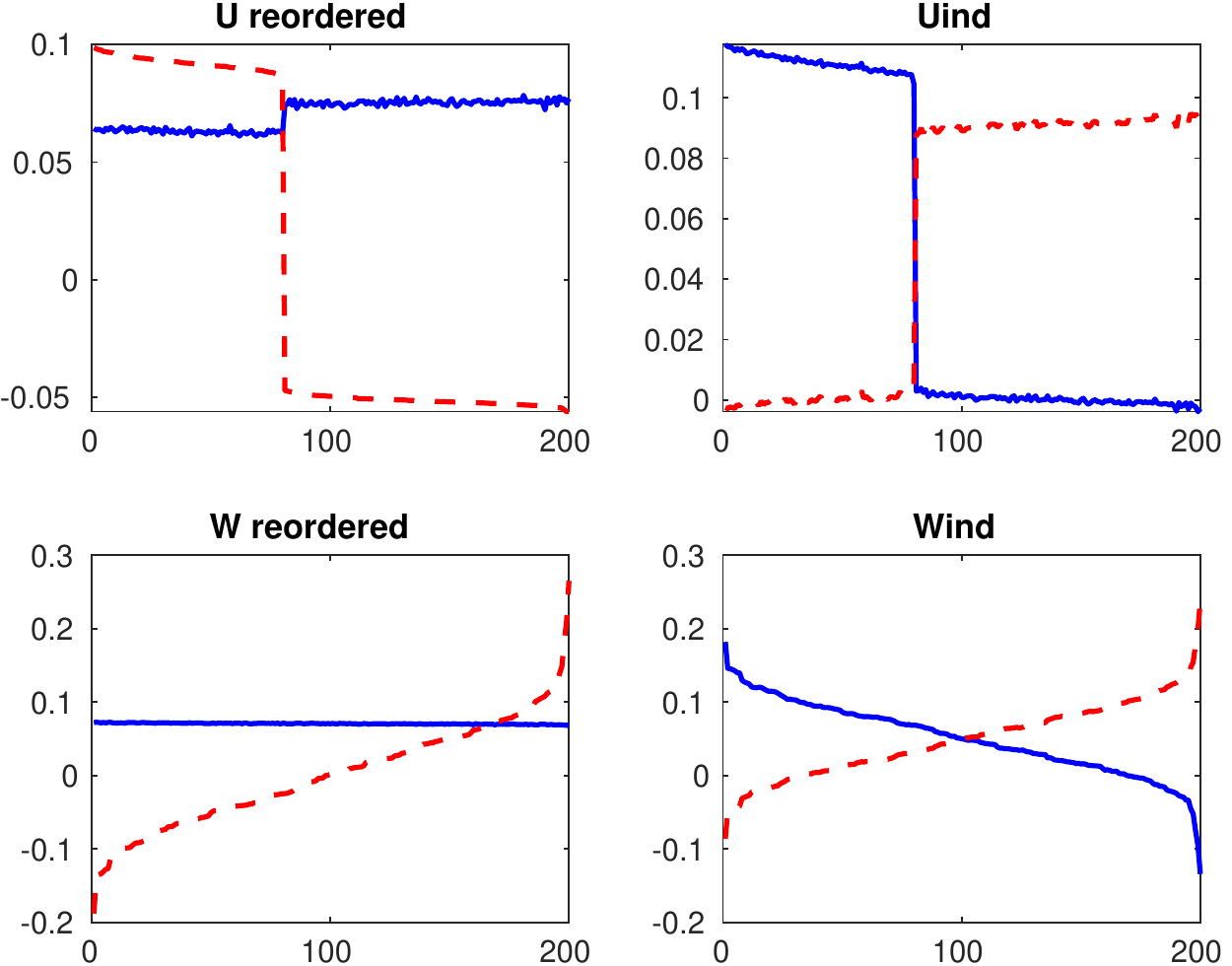}
\caption{Example \ref{ex:Ex1}, rank-(2,2,2) approximation. The column
  vectors of the solution 
   $U$ after reordering (upper left), 
   indicator 
  vectors for $U$ (upper right),  reordered $W$ vectors (lower left),
  and computed indicator  vectors for 
  $W$ (bottom right).} 
\label{fig:VectorsEx1-222}
\end{figure}
Thus we can safely conclude that the correct structure was  identified
using the rank-(2,2,1) approximation.
\end{example}
 
\medskip

\begin{example}\label{ex:Ex2}
  The tensor was constructed with  structure close to
  that of
  Proposition \ref{prop:0BB0}. The results after computing a
  rank-(2,2,1) approximation are given in Figures
  \ref{fig:VectorsEx2}-\ref{fig:ReconTensorEx2}. 

  %
%
  %
\begin{figure}[htbp!]    
\centering
\includegraphics[width=.65\textwidth]{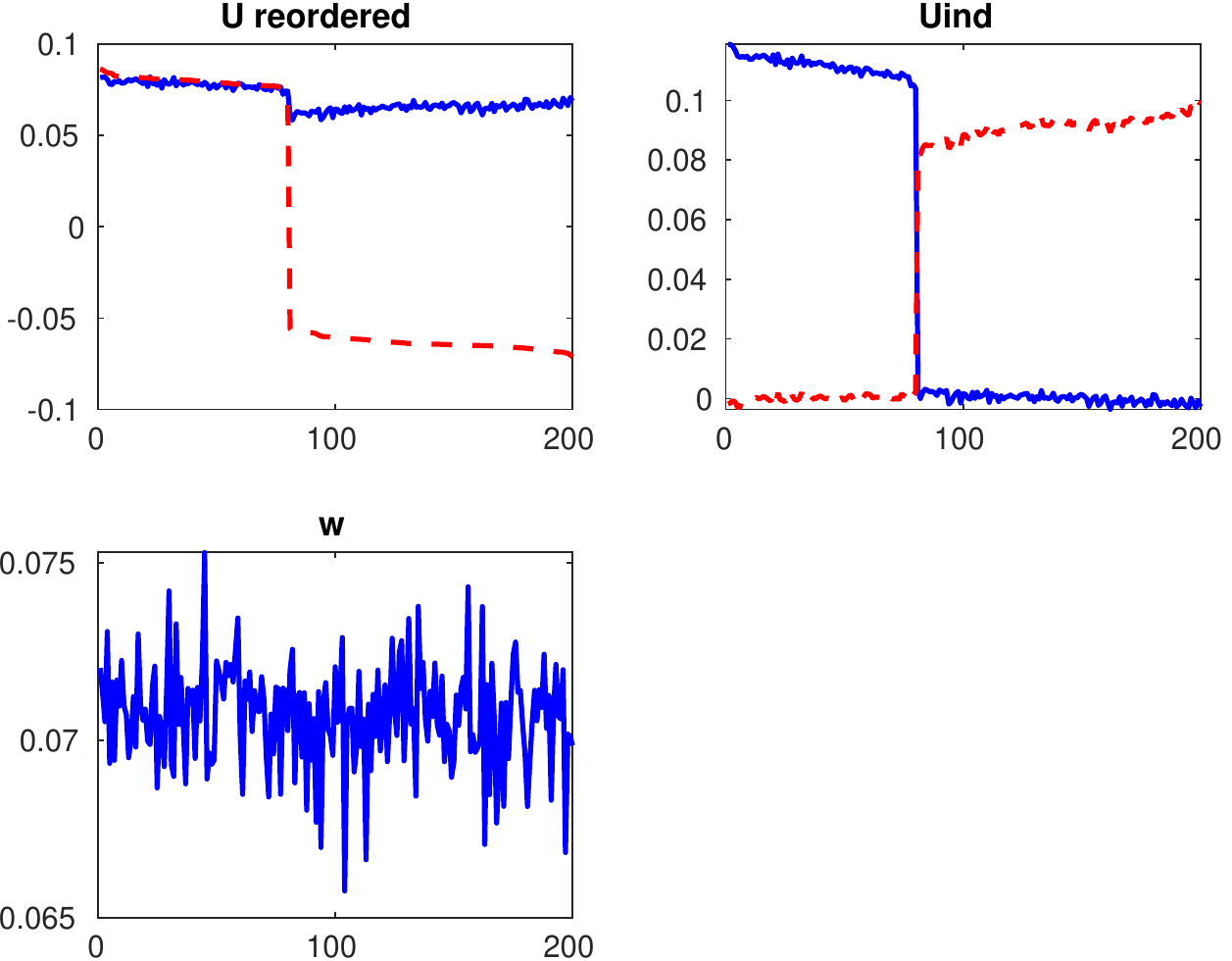}
\caption{Example \ref{ex:Ex2}, rank-(2,2,1) approximation. The
  reordered column vectors of 
   $U$ (upper left), indicator vectors for $U$ 
   (upper right),  $w$ (lower left).} 
\label{fig:VectorsEx2}
\end{figure}

\begin{figure}[htbp!]    
\centering
\includegraphics[width=.6\textwidth]{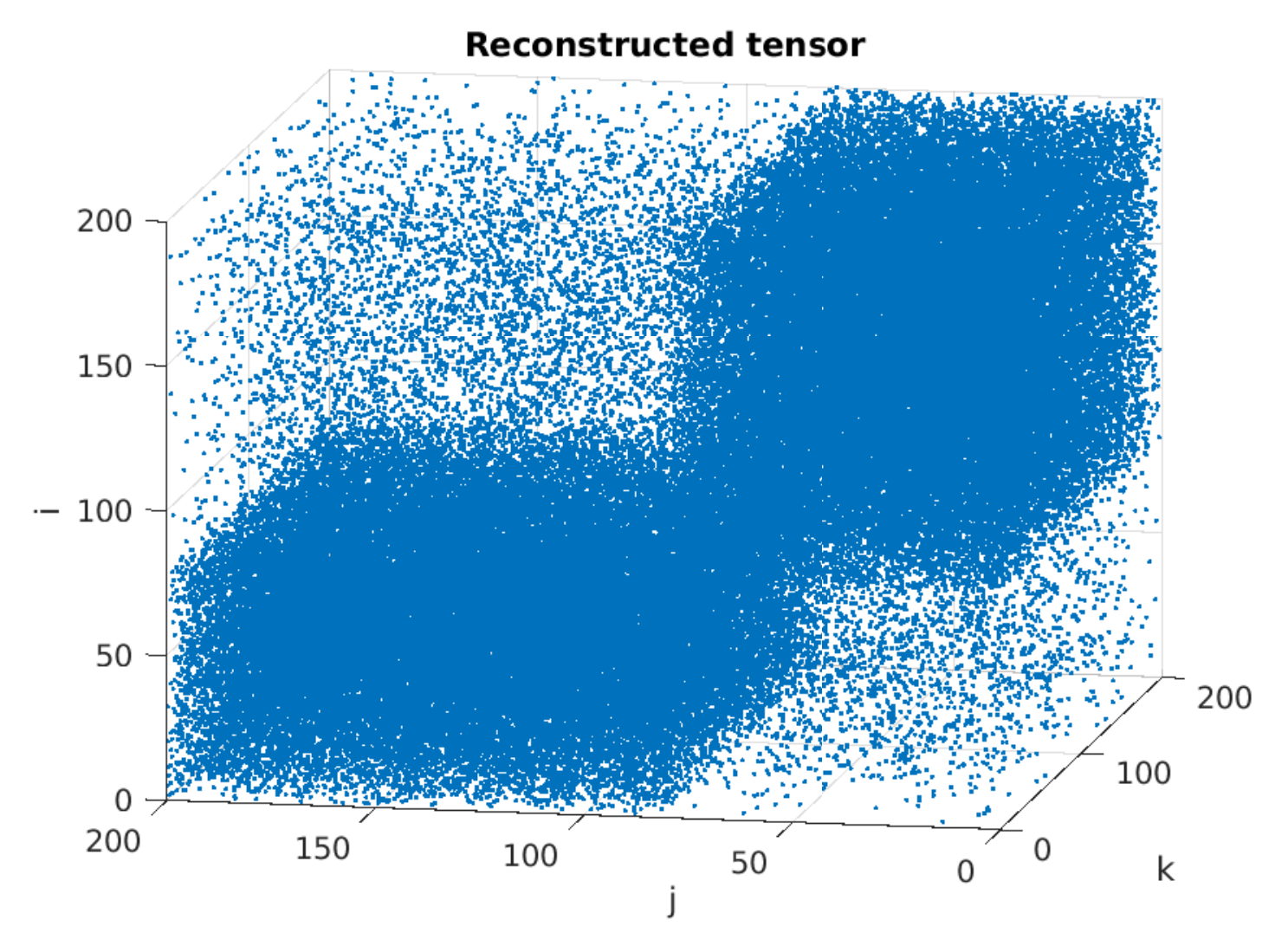}
\caption{Example \ref{ex:Ex2}. The reconstructed  tensor.}
\label{fig:ReconTensorEx2}
\end{figure}
The  structure tensor was
\[
    \cF_{\mathrm{struct}} =
    \begin{pmatrix} 
-0.62 &   11.46 \\
	   11.46 &  -1.16
    \end{pmatrix}, 
  \]
  which is a clear indication for the structure of Proposition
  \ref{prop:0BB0}. The  norm of the
  original tensor $\cA$ was 111.13 and the norm of $\cF$ was 16.34;
  the norm of the core 
  tensor in a rank-(2,2,2) approximation was exactly the same. This
  shows that there was no significant mode-3 structure that was missed
  in the rank-(2,2,1) approximation. 
\end{example}

\medskip

\begin{example}\label{ex:Ex3}
  We constructed a tensor with  structure close to that  in
  \eqref{eq:A-1,2-red-1}-\eqref{eq:A-1,2-red-2}. Rank-(2,2,1) and
  rank-(2,2,2) approximations were computed with $\| \cF \| = 15.03$
  and $\| \cF \| = 15.96$, respectively. Clearly there was significant
  mode-3 structure that was visible in the rank-(2,2,2) approximation.
  However, with rank-(2,2,3) the norm of the core tensor did not
  become larger.

  We
  reordered the vectors and applied the reorderings to the tensor. For
  visibility reason we do not give a spy plot of the tensor but of a
  slice, see Figure \ref{fig:ReconSliceEx4}. 
\begin{figure}[htbp!]    
\centering
\includegraphics[width=.4\textwidth]{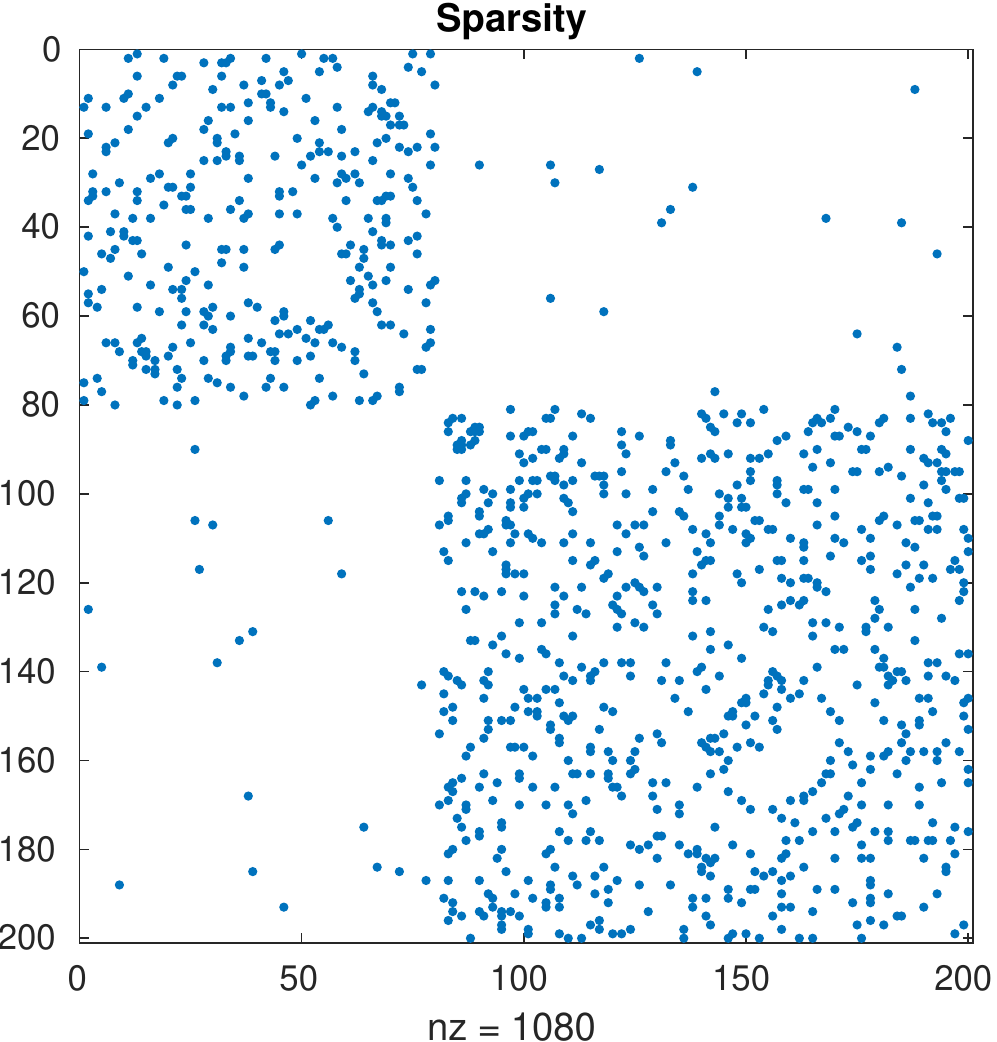}
\includegraphics[width=.4\textwidth]{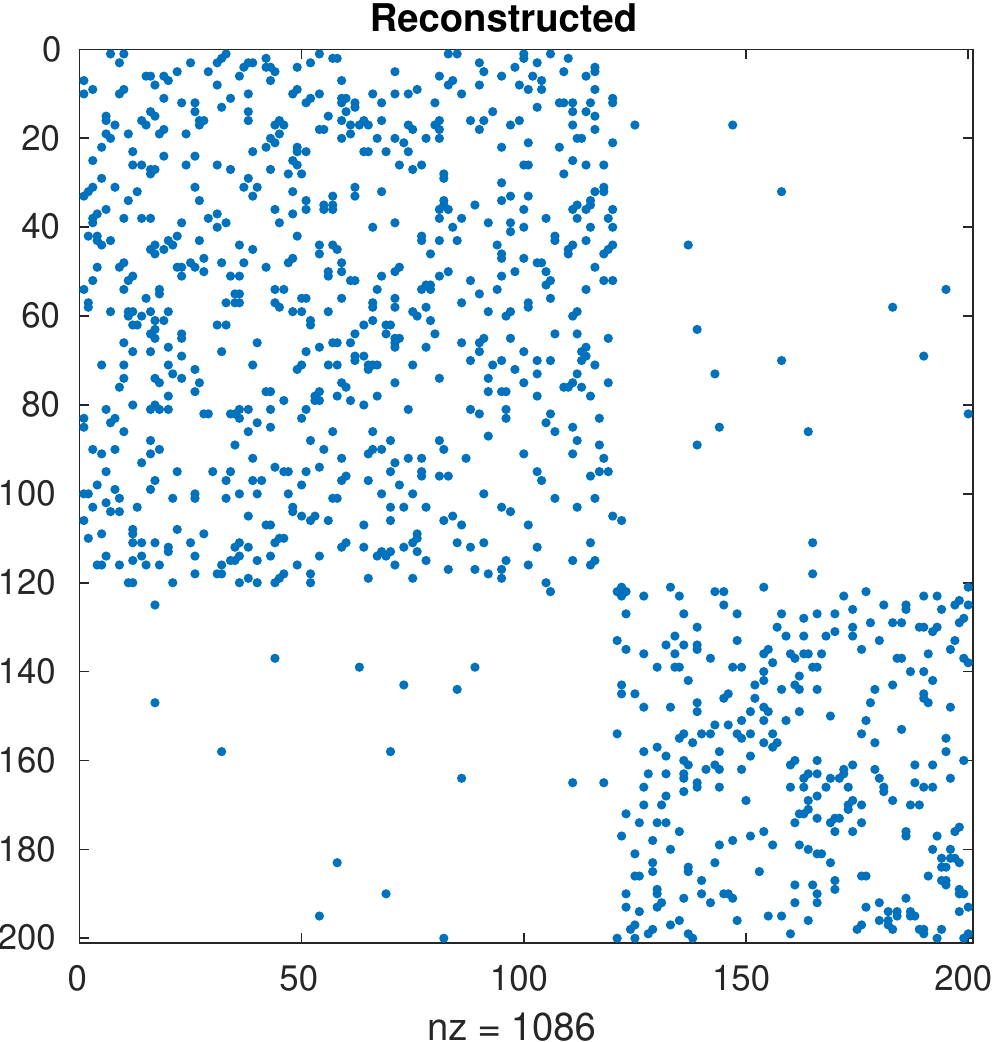}
\caption{Example \ref{ex:Ex3}. Slices of the original tensor (left)
  and the reconstructed  tensor (right).}
\label{fig:ReconSliceEx4}
\end{figure}
The vectors are illustrated in Figure \ref{fig:VectorsEx4}. 
\begin{figure}[htbp!]    
\centering
\includegraphics[width=.65\textwidth]{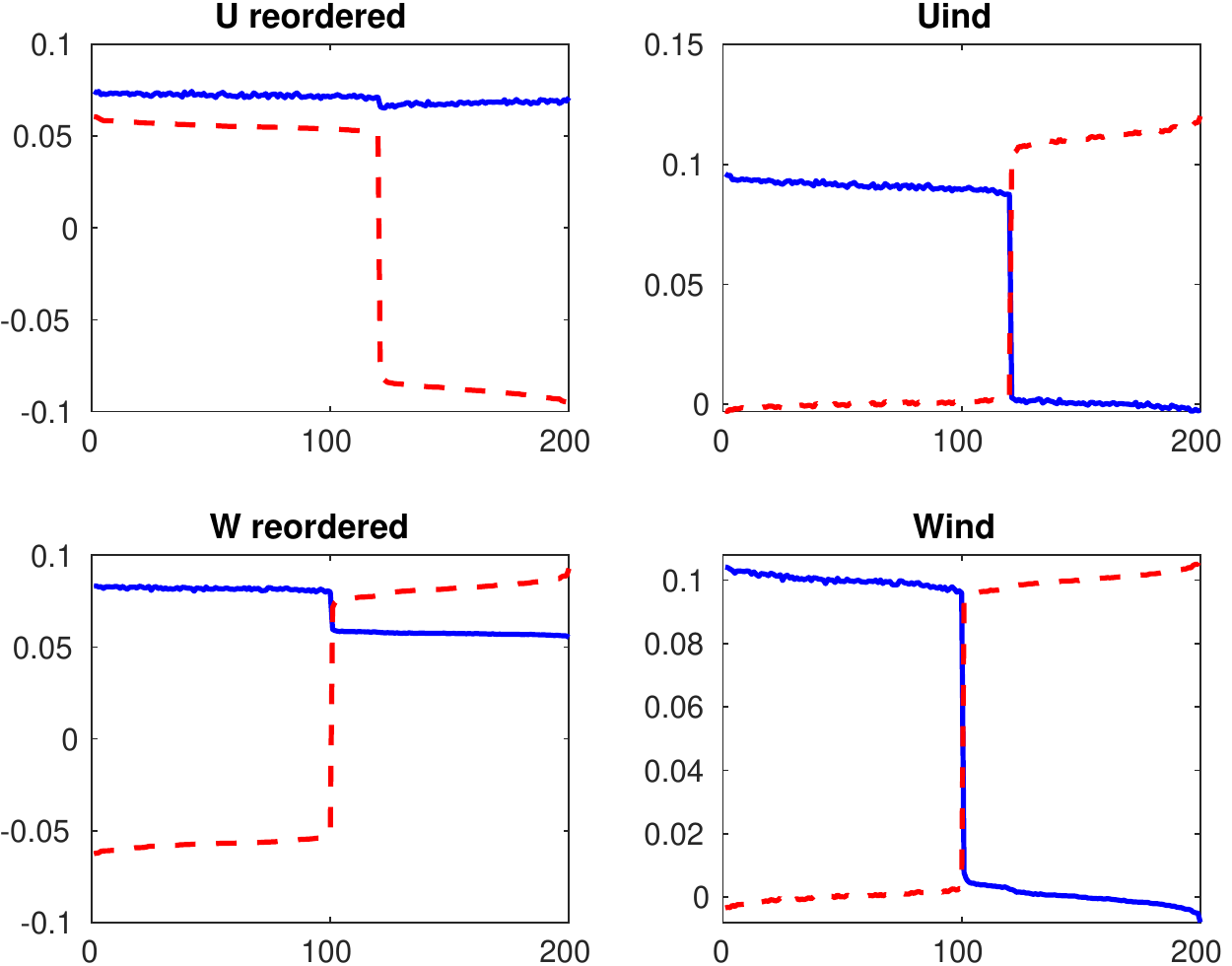}
\caption{Example \ref{ex:Ex3}, rank-(2,2,2) approximation. The
  reordered column vectors of 
   $U$ (upper left), indicator vectors for $U$ 
   (upper right),  reordered vectors of $W$ (lower left), and
   indicator vectors for $W$ (lower right).} 
\label{fig:VectorsEx4}
\end{figure}
The structure tensor was 
\begin{align*}
  \cF_{\mathrm{struct}}(:,:,1)&=
  \begin{pmatrix}
      9.33  &  -0.36\\
	    -0.36 &    8.83
  \end{pmatrix}, \quad 
  \cF_{\mathrm{struct}}(:,:,2) =
  \begin{pmatrix}
	    5.92 &    -4.86\\
	    -4.86 &   3.91
  \end{pmatrix}.
\end{align*}
We also ran a similar example, where the tensor was exactly
reducible. The structure tensor was
\[
    \cF_{\mathrm{struct}}(:,:,1)=
  \begin{pmatrix}
      13.31  &  -0.031\\
	    -0.031 &    7.79
  \end{pmatrix}, \quad 
  \cF_{\mathrm{struct}}(:,:,2) =
  \begin{pmatrix}
	    3.04 &    -0.31\\
	    -0.31 &   -5.21
  \end{pmatrix}.
\]
The singular values of the matrix $\widehat U$ in \eqref{eq:Uhat} were
$1.414$ and $0.010$, and those of the corresponding $\widehat W$  were
$1.414$ and 0.025. This shows clearly that for the exactly reducible
tensor the indicator vectors are only in approximate indicator
form. However, they give the correct partitioning, also in the almost
reducible case. 
\end{example}
\bigskip

\begin{example}\label{ex:Ex4}
  A  tensor with the structure close to \eqref{eq:reducible-tensor3} 
  (i.e., almost 3-reducible) was constructed and analyzed. The results
  were completely analogous   to those in Example \ref{ex:Ex3}, and we
  therefore omit them here. 
\end{example}

\subsection{Real Data}
\label{sec:real-data}
As we remarked earlier it is probably  rare that tensor data are
close to reducible. However, best low-rank approximation can be used 
for analysis by clustering. This is the topic of \cite{eldehg20c}. We
here give a brief account of an example from that paper. 

\medskip

\begin{example}\textbf{Clustering of news texts}. \label{ex:Reuters}
  The tensor is a collection of graphs constructed from a corpus of
  news texts published by Reuters during 66 days after the September
  11 attack on the U.S. The vertices are words; there is an edge
  between two words if and only if they appear in the same text unit
  (sentence). The weight of an edge is its frequency. The tensor has
  dimension $13332 \times 13332 \times 66$ and is
  (1,2)-symmetric. Naturally it is  impossible 
  to discern any structure by visual inspection of a spy plot of the tensor.

  We computed  a best rank-(2,2,2) approximation of the tensor using
  the BKS method \cite{eldehg20b}, where the outer iterations were
  intialized with 10 HOOI iterations for rank-(1,1,1). With
  relative termination criterion $10^{-8}$,  30
  iterations were required, and the execution time was about 44
  seconds.
  The
  reordered $U$ and $W$ vectors are shown in Figure
  \ref{fig:UW66reordered}.
\begin{figure}[htbp!]    
\centering
\includegraphics[width=.6\textwidth]{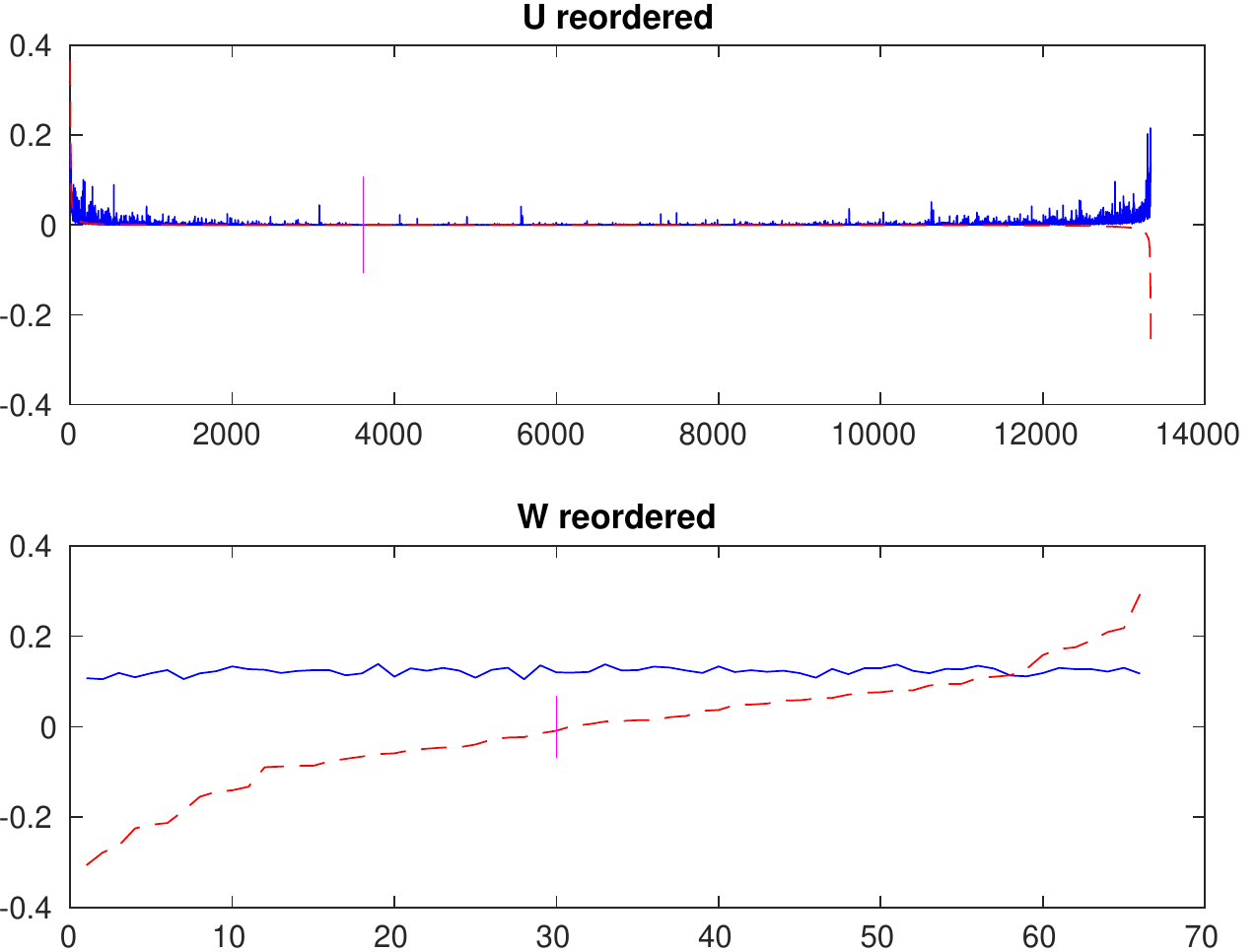}
\caption{ The column vectors of $U$ and $W$ after
  reordering; $u_1$ and $w_1$ are blue and the others red. The magenta
  vertical lines show  the position, where the elements of the reordered
  vectors $u_2$  and $w_2$   change sign.  }  
\label{fig:UW66reordered}
\end{figure}
After the reordering of the vectors (words), the two main topics of
the texts are at the beginning and at the end, see Table
\ref{tab:terms66}. 
\begin{table}[htbp!]
  \centering
  \caption{Beginning terms, middle terms and end terms after
    reordering. The table is organized so that the most significant
    terms appear at the top of the first column and the bottom of the
    third (consistent with the reordered vectors). The middle column
    shows a few terms 
    at the middle, which are relatively insignificant with respect to
    the topics in the other two columns.\label{tab:terms66} }
  \begin{tabular}{ll|ll|llll}
    &Topic 1& &Insignificant & &Topic 2&\\
    \hline
     &world\_trade\_ctr&    &planaria&         &hand&  \\      
    &pentagon&           &mars&             &ruler&       \\
    &new\_york&           &i.d&              &saudi&     \\  
    &attack&             &eyebrow&          &guest&    \\   
    &hijack&             &calendar&         &mullah&   \\ 
    &plane&              &belly&            &movement&   \\ 
    &airliner&           &auction&          &camp&   \\
    &tower&              &attentiveness&    &harbor& \\   
    &twin&               &astronaut&        &shelter&  \\   
    &washington&         &ant&              &group&      \\ 
    &sept&               &alligator&        &organization&\\
    &suicide&            &adjustment&       &leader&      \\
    &pennsylvania&       &cnd&              &rule&   \\
    &people&             &nazarbayev&       &exile&    \\   
    &passenger&          &margrit&          &guerrilla&  \\ 
    &jet&                &bailes&           &afghanistan& \\
    &mayor&              &newscast&         &fugitive&    \\
    &110-story&          &turgan-tiube&     &dissident&   \\
    &mayor\_giuliani&     &N52000&           &taliban&     \\
    &hijacker&           &ishaq&            &islamic&     \\
    &commercial&         &sauce&            &network&     \\
    &aircraft&           &sequential&       &al\_quaeda&   \\
    &assault&            &boehlert&         &militant&    \\
    &airplane&           &tomato&           &saudi-born&  \\
    &miss&               &resound&          &bin\_laden&    
  \end{tabular}
\end{table}
A large number of insignificant words are at the middle (corresponding
to small components of the vector $u_1$). Notice that these words can
not be considered as stop words in the sense of information retrieval
\cite{byrn:99,bebr:05}, but they are certainly insignificant with
respect to the task of finding the main topics in the texts.

We divided the 66 days into two groups, where the first, $D_A$ (30
days), corresponds to the (original) indices of $w_2$ to the left of
the sign change, see Figure \ref{fig:UW66reordered}. The second group,
$D_B$ (36 days), corresponds to the the (original) indices to the
right.  It turns out that the reordered days are not consecutive:
\medskip\medskip 
{\small
\begin{verbatim}
   DA: XXXXXXXXXX X  XX X  XX XX    X XX           XXX  X X X X   X X    
   DB:           X XX  X XX  X  XXXX X  XXXXXXXXXXX   XX X X X XXX X XXXX 
\end{verbatim}
}
\medskip\medskip

The reordered tensor is shown in Figure \ref{fig:reuters66reordered}. 
\begin{figure}[htbp!]    
\centering
\includegraphics[width=.6\textwidth]{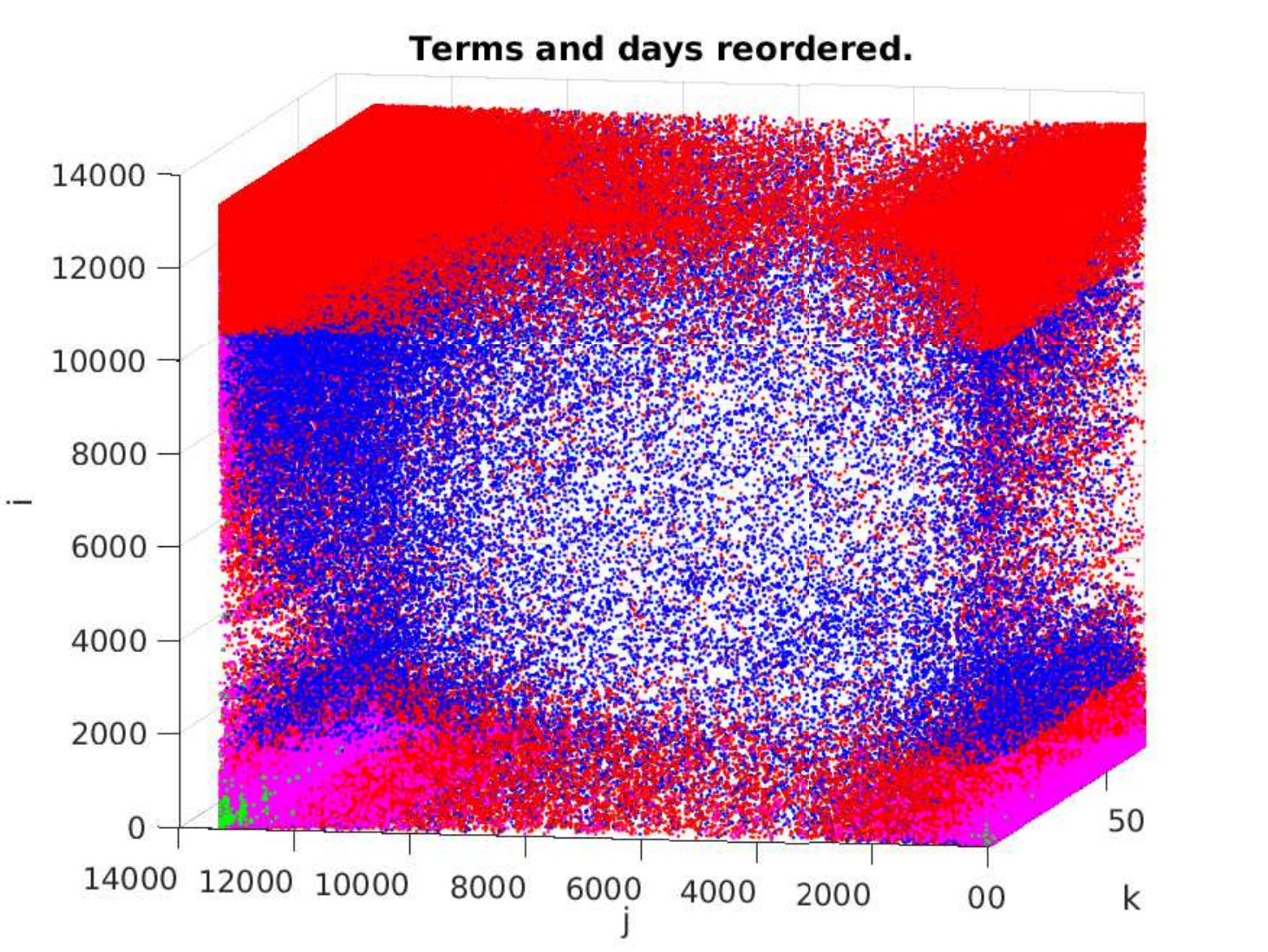}
\caption{The reordered Reuters tensor.  }  
\label{fig:reuters66reordered}
\end{figure}
It is seen that the ``mass'' is concentrated to the four corners in
term modes, which correspond to Topics 1 and 2 in Table
\ref{tab:terms66}. We computed the norms of the subtensors of
dimensions $1000 \x 1000 \x 66$ at the corners (cf. the shape of $U$ in
Figure \ref{fig:UW66reordered}), and they were $50.6$, $25.5$, $25.5$,
and $35.1$, respectively. As the norm of $\cA$ was $128.6$, about 55\%
of the total mass was concentrated at  those corners. This shows that
the terms of Topic 1 cooccurred to a large extent with terms within
the same topic, and similarly for Topic 2. In addition, there was
considerable cooccurrence between terms from the two main topics.

Clearly, the tensor of this example is far from reducible.
Still, the example indicates that the partitioning given by the
rank-$(2,2,2)$ approximation orders the vertices of the cooccurrence
graphs so that the dense parts of the graphs are at both ends.
The middle part is sparser, corresponding to terms that can be
considered as noise (in \cite{eldehg20c} we removed the 9000 middle
terms, and repeated the analysis; the same topics were produced).  
In \cite{eldehg20c} we partitioned the reordered tensor, computed new
rank-(2,2,2) approximations of subtensors, and found subtopics of those in Table
\ref{tab:terms66}. 
\end{example}

\section{Conclusions}
\label{sec:conclusions}
In this paper we have presented a generalization of spectral graph
partitioning to sparse tensors, where the 3-slices are
normalized adjacency matrices of undirected graphs. We have shown that
for four cases of tensors that are close to reducible, it is possible to
determine the reducibility structure by computing a low-rank best
multilinear approximation of the tensor.  
%
%
The algorithm is applied to synthetic data and it is shown that the
reducibility structure is accurately computed.

Due to the presence of noise, it probably rare to find tensors in
applications that are close to reducible. We give one example with a
tensor from a text analysis application, where actually most of the
tensor data are noise, and there is no sharp distinction between
useful information and noise.  Here the best low-rank approximation
performs a clustering that separates
the relevant information from the noise.  This property is likely to
make the methods of this paper useful in applications. In
\cite{eldehg20c} we investigate the use of the methods in a few
different data science applications. It is shown that the method can
be used to extract information from large tensors with noisy data.

\bigskip

\appendix

\section{Nonnegativity of the $\mathbf{Rank}$-(1,1,1) Solution for
  Nonnegative  
  Tensors}
\label{app:111}

Using Perron-Frobenius theory for nonnegative matrices it can be
proved that the eigenvector corresponding to the dominant eigenvalue
is nonnegative \cite[Chapter 8]{mey:00}. Generalizations to tensors
are given in \cite{fgh13,gahe16}, where the norms are $p$-norms with
restrictions on the value of $p$\footnote{Note that the counterexamples in
\cite{fgh13} do not imply the nonexistence of a nonnegative solution
of the best $\rank$-(1,1,1) approximation problem for a nonnegative
tensor.}.  Even if the result seems to be well-known, we are
not aware of any proof in the literature that the best rank-(1,1,1)
approximation of a nonnegative tensor is nonnegative. 

\begin{lemma}\label{lem:nonneg}
  Let $\cA \in \RR^{l \times m \times n}$ be a nonnegative tensor. Then 
  the best rank-(1,1,1) approximation problem has a nonnegative
  solution.
  \end{lemma}

  \begin{proof}
         Let
    $(u,v,w)$ be  a solution. Without loss of
    generality we can assume that
    \[
      \max_{\| x \| = \| y \| = \| z \|=1} | \tmr{\cA}{x,y,z} | =
      | \tmr{\cA}{u,v,w} | = 1, 
    \]
and  that
  \[
    u = 
    \begin{pmatrix}
      u_1 \\ u_2
    \end{pmatrix}, \qquad
    v =
    \begin{pmatrix}
      v_1 \\ v_2 
    \end{pmatrix}, \qquad
    w =
    \begin{pmatrix}
      w_1 \\      w_2
    \end{pmatrix},
  \]
  where $u_1 > 0$, $v_1 > 0$, $w_1 > 0,$ and $u_2 \leq 0$,
  $v_2 \leq 0$, and $w_2 \leq 0$.  Assume that the tensor is
  partitioned accordingly in the three modes and denote the subtensors
  $\cA_{\lambda\mu\nu}$, for $\lambda,\mu,\nu=1,2$. Let $| x |$ denote the vector of absolute
  values of the components of $x$, and  let $u \otimes v \otimes w \in
  \RR^{l \times m \times n}$ denote the tensor (outer) 
  product of the vectors\footnote{Using our notation for tensor-matrix
    multiplication, we have
    $u \otimes v \otimes w = \tml{u,v,w}{\mathbbm{1}}$, where
    $\RR^{1 \times 1 \times 1} \ni \mathbbm{1}$ is a tensor.}. Then
    \begin{align}\label{eq:cA1}
      \| \cA - u \otimes v \otimes w \|^2 
      &= \| \cA_{111}-u_1 \otimes  v_1 \otimes w_1 \|^2
      + \| \cA_{121} + u_1 \otimes | v_2 | \otimes w_1 \|^2 \\
      &+ \qquad \cdots \nonumber\\
      &+ \| \cA_{212} - | u_2 | \otimes  v_1   \otimes | w_2 | \|^2
        + \| \cA_{222} + | u_2 | \otimes | v_2 | \otimes | w_2 | \|^2,
        \label{eq:cA2}
    \end{align}
    where we have omitted four similar terms. Since all subtensors
    are nonnegative, we have, e.g.,
    \begin{equation}
      \label{eq:A121-uvw}
           \| \cA_{121} + u_1 \otimes | v_2 | \otimes w_1 \|^2 \geq
      \| \cA_{121} - u_1 \otimes | v_2 | \otimes w_1 \|^2, 
    \end{equation}
    The only situation when strict inequality can occur is if a
    component of $\cA_{121}$ is positive, and the corresponding
    component of $v_2$ is nonzero; however, that would be a
    contradiction to the assumption that $(u,v,w)$ is a
    solution. Thus, we have equality in \eqref{eq:A121-uvw}. It
    follows that we can replace all $+$ in
    \eqref{eq:cA1}-\eqref{eq:cA2} by $-$, and have
    \begin{equation}
      \label{eq:uvw-ineq}
      \| \cA - u \otimes v \otimes w \|^2  =
      \| \cA - x \otimes y \otimes z \|^2,      
    \end{equation}
    where
    \[
          x = 
    \begin{pmatrix}
      u_1 \\ | u_2 |
    \end{pmatrix}, \qquad
    y =
    \begin{pmatrix}
      v_1 \\ | v_2 |
    \end{pmatrix}, \qquad
    z =
    \begin{pmatrix}
      w_1 \\ | w_2 |
    \end{pmatrix}.
  \]
Thus we have a nonnegative solution $(x,y,z)$. 
  \end{proof}

  \medskip
  

\section{Lemma for a Rank-2 Matrix}
\label{app:C12}

\begin{lemma}\label{lem:rank-2}
  Let $U \in \RR^{m \times 2}, \; m > 2$, satisfying $U\tp U = I$, and
  let   $A \in \RR^{2 \times 2}$ be symmetric. Put
    \[
      Z = U A U\tp =
      \begin{pmatrix}
        Z_{11} & Z_{12}\\
         Z_{12} & Z_{22}
       \end{pmatrix}, \qquad
       Z_{12} \in \RR^{m_1 \times m_2}, \quad m_1 + m_2 =m,
     \]
     and let the columns of $U$ be partitioned accordingly.
     
     Then $ Z_{12}=0$ and rank$(Z)=2$  if and only if 
     \begin{equation}
       \label{eq:Z-gen}
       U=
       \begin{pmatrix}
         u_1 & 0\\
         0   & u_2
       \end{pmatrix}, \quad 
       A = \begin{pmatrix}
         a_{11} & 0 \\
         0   & a_{22}
       \end{pmatrix}, 
       \quad
       \rank(A)=2, 
     \end{equation}
     or 
     \begin{equation}\label{eq:Z-12}
       U=
       \begin{pmatrix}
         u_{11} & u_{12}\\
         0 & 0
       \end{pmatrix}, \quad
       \rank(A)=2, 
       \end{equation}
   or
       \begin{equation}\label{eq:Z-22}
       U=
       \begin{pmatrix}
         0 & 0 \\
         u_{21} & u_{22}
       \end{pmatrix}, \quad         \rank(A)=2.
      \end{equation}
\end{lemma}

\begin{proof}
  First note that, for any orthogonal $Q
  \in \RR^{2 \times 2}$,
  \begin{equation}
    \label{eq:Z-Grassmann}
    Z = (UQ) (Q\tp A Q) (UQ)\tp,
  \end{equation}
  so the factors   in the product are not   unique. Let
  \[
    A=
    \begin{pmatrix}
      a_{11} & a_{12}\\
      a_{12} & a_{22}
    \end{pmatrix},
    \qquad
    U=
    \begin{pmatrix}
      x_1 & y_1\\
      x_2 & y_2
    \end{pmatrix},
  \]
  with the obvious partitioning of the column vectors in $U$. 
  By straightforward computation
  \begin{equation}
    \label{eq:C-12}
    Z_{12} = a_{11} x_1 x_2\tp + a_{22} y_1 y_2\tp + a_{12}(y_1 x_2\tp + x_1
    y_2\tp).
  \end{equation}
  Thus, any of the U-structures in \eqref{eq:Z-gen} or \eqref{eq:Z-12}
  implies $Z_{12}=0.$

  Then assume that $Z_{12}=0$, $x_1 \neq 0$ and $x_2 \neq 0$. By
  multiplying \eqref{eq:C-12} from the 
  left and right by $x_1\tp$ and $x_2$, respectively, we see that
  \[
    y_1 = \tau_1 x_1, \quad y_2 = \tau_2 x_2,
  \]
  for some $\tau_1$ and $\tau_2$. The requirements that the column
  vectors are orthogonal and have length 1, give
  \[
    U=
    \begin{pmatrix}
      x_1 & \tau x_1\\
      x_2 & -\frac{1}{\tau} x_2
    \end{pmatrix},
    \qquad
    \tau=\frac{\|x_2\|}{\|x_1\|}.
    \] 
    Inserting this in \eqref{eq:C-12} we get
    \[
      0=Z_{12} = \gamma  \, x_1 x_2\tp:= (a_{11} - a_{22} +
      a_{12}(\tau -\frac{1}{\tau})) \,  x_1 x_2\tp,
    \]
    which implies $\gamma=0.$ 
    So we have
    \[
      Z = U A U\tp =
      \begin{pmatrix}
        x_1 & \tau x_1 \\
        x_2 & - \frac{1}{\tau} x_2
      \end{pmatrix}
      \begin{pmatrix}
        a_{11} & a_{12}\\
        a_{12} & a_{22}
      \end{pmatrix}
      \begin{pmatrix}
        x_1\tp & x_2\tp\\
        \tau x_1\tp & - \frac{1}{\tau} x_2\tp 
      \end{pmatrix}. 
    \]
    We now insert the orthogonal matrix $Q$,
    \[
      Z = (UQ) (Q\tp A Q) (UQ)\tp, \qquad
      Q=
      \begin{pmatrix}
        \| x_1 \| & \| x_2 \| \\
        \| x_2 \| & - \| x_1 \| \\
      \end{pmatrix}, 
    \]
    and put $\widehat{U}=UQ$ and $ \widehat{A} = Q\tp A Q$. 
    Straightforward computation gives 
    \[
      \widehat{U} = 
      \begin{pmatrix}
        \frac{1}{\| x_1 \|} \,  x_1 & 0 \\
        0 & \frac{1}{\| x_2 \|} \, x_2
      \end{pmatrix}, \quad
      \widehat A= 
      \begin{pmatrix}
        \hat a_{11} & \hat a_{12}\\
        \hat a_{12} & \hat a_{22}
      \end{pmatrix},
    \]
    where $\hat a_{12}=-\| x_1 \| \, \| x_2 \| \, (a_{11} - a_{22} +
    a_{12}(-\frac{1}{\tau} + \tau)) = 0$. This means that $U$ and $A$ 
    have the structure  \eqref{eq:Z-gen}.

    Assume that $Z_{12}=0$ and $\rank(Z)=2,$  $x_1 \neq 0$ and $y_1 \neq 0$, and
    $x_2=y_2= 0$. Then it is easy to show that \eqref{eq:Z-12}
    holds. Analogously one can easily prove the case of  \eqref{eq:Z-22}.

    Finally we must check the following case: assume that $Z_{12}=0$,
    $x_1 \neq 0$ and $y_1 \neq 0$, 
    $x_2=0$,and $y_2 \neq 0$. From orthogonality of the columns of $U$
    we have $x_1\tp y_1 =0$, and using
    \[
      0= Z_{12}= (a_{22} y_1 + a_{12} x_1) y_2\tp, 
    \]
    we get     $a_{22}=a_{12}=0$. Therefore, in this case $Z$ is 
    independent of $y_1$ and   $y_2$, and 
    \[
      Z=
      \begin{pmatrix}
        a_{11} x_1 x_1\tp & 0 \\
        0      & 0
      \end{pmatrix},
    \]
    which means that rank$(Z)=1$. The same is true three  similar cases.   
\end{proof} 

\medskip

Note that the case when
\[
U =
\begin{pmatrix}
  0 & u_{12}\\
  u_{21} & 0
\end{pmatrix},
\]
is equivalent to \eqref{eq:Z-gen} due to the property
\eqref{eq:Z-Grassmann}.

In the case \eqref{eq:Z-gen} the blocks $Z_{11}$ and $Z_{22}$ both
have rank 1.  In the cases \eqref{eq:Z-12} and \eqref{eq:Z-22}, 
$Z_{11}$  and   $Z_{22}$ have rank 2, respectively.

\medskip

\bibliography{/media/lars/ExtHard/WORK/forskning/BIBLIOGRAPHIES/general,/media/lars/ExtHard/WORK/forskning/BIBLIOGRAPHIES/LE-papers}
\bibliographystyle{plain}

\end{document}